\newcommand\modulename{PlainPaper}
\newcommand\papername{Compute}
\newcommand\papertitle{A computational proof of the linear Arithmetic Fundamental Lemma of $\GL_4$}
\def\href#1{}
\documentclass[12pt]{amsart}

\usepackage{pdfsync}
\usepackage{fullpage}
\usepackage{python}
\usepackage{mathrsfs}
\usepackage{latexsym, amssymb, amsmath, amscd, amsthm, amsxtra}
\usepackage{eucal}
\usepackage[utf8]{inputenc}
\usepackage{verbatim}
\usepackage{hyperref}
\usepackage[english]{babel}
\usepackage{diagbox}
\usepackage{mathtools}
\usepackage{todonotes} 
\usepackage{cite}
\usepackage{graphicx}
\usepackage{amssymb}
\usepackage{epstopdf}
\usepackage{tikz}
\usetikzlibrary{calc}
\usetikzlibrary{graphs}
\usepackage[all]{xy}
\usepackage{color}
\usepackage{fancyhdr}
\usepackage{arydshln}
\usepackage{times}
\usepackage{datetime}
\usepackage{scrtime}
\usepackage{CJKutf8}
\usepackage[refpage]{nomencl}
\usepackage{cancel}
\usepackage{ulem}
\usepackage{bbm}

\def\load#1{\input{\modulename/Command/#1}}
\def\Load#1{\input{\modulename/Command/#1}}
\long\def\[#1]#2{ \begin{#1}#2\end{#1}}

\renewcommand{\[}{\left[}

\newcommand{\id}{\mathrm{id}}

\newcommand{\ZZ}{\mathbb{Z}}

\newcommand{\PP}{\mathbb{P}}

\newcommand{\map}[3]{#1:#2\longrightarrow #3}
\newcommand{\mapp}[2]{#1\longrightarrow #2}

\newcommand{\maps}[5]{\begin{array}{llrll}
#1&:&#2&\longrightarrow &#3\\\\
&&#4&\longmapsto &#5\\
\end{array}}
\newcommand{\mapi}[4]
{\begin{array}{lrrll}
&#1&\longrightarrow &#2\\\\
&#3&\longmapsto &#4\\
\end{array}}
\newcommand{\mapa}[3]
{\xymatrix{{#2}\ar[rr]^{#1}&&{#3}}}

\DeclareMathOperator{\Ker}{\mathrm{Ker}}

\newcommand{\CC}{\mathbb{C}}

\newcommand{\dd}{\mathrm{d}}

\newcommand{\RR}{\mathbb{R}}
\newcommand{\HH}{\mathbb{H}}

\newcommand{\CFF}[1]{\overline{\FF}_#1}
\newcommand{\FF}{\mathbb{F}}

\newcommand{\one}{\mathbbm{1}}
\newcommand{\OO}[1]{\mathcal{O}_{#1}}

\newcommand{\vv}{\textbf{v}}

\newcommand{\ep}{\boldsymbol{\epsilon}}

\DeclareMathOperator{\End}{\mathrm{End}}
\newcommand{\Vol}{\mathrm{Vol}}
\DeclareMathOperator{\Hom}{\mathrm{Hom}}

\DeclareMathOperator{\Spf}{\mathrm{Spf}\text{ }}

\newcommand{\Stab}{\mathrm{Stab}}

\newcommand{\Isom}{\mathrm{Isom}}

\newcommand{\Isog}{\mathrm{Isog}}

\DeclareMathOperator{\length}{\mathrm{length}}
\DeclareMathOperator{\Aut}{\mathrm{Aut}}
\DeclareMathOperator{\gl}{\mathrm{Mat}}
\newcommand{\gr}{\mathrm{Gr}}
\DeclareMathOperator{\GL}{\mathrm{GL}}

\DeclareMathOperator{\im}{\mathrm{Im}}

\newcommand{\ceil}[1]{\left\lceil #1\right\rceil}

\newcommand{\ezer}[4]{\xymatrix{
{#1}\ar@<-.5ex>[rr]_{#3}\ar@<.5ex>[rr]^{#4}&&{#2}}
}

\input{\papername/linear}

\long\def\[#1]#2{ \begin{#1}#2\end{#1}}

\begin{document}
\def\load#1{\input{\modulename/Style/#1}}
\def\Load#1{\input{\modulename/Style/#1}}
\let\yourlabel=\label
\def\lb#1{\yourlabel{#1}\comments{{\color{blue}\textbf{Lable}}:\textbf{#1   }  }}
\def\status{\use{Status}}
\def\test{t}
\def\comments#1{\ifthenelse{\equal\status\test}
{{\color{red} #1}}
{}
}

\long\def\s#1{\ifthenelse{\equal\status\test}
{{\color{blue} \textbf{{\color{red}]}#1 \color{red}[}}}
{}
}

\long\def\ss#1{\ifthenelse{\equal\status\test}
{{\color{purple} \textbf{{\color{red}\}}#1 \color{red}\{~\\}}}
{}
}

\newcommand{\commenta}[1]{{\color{blue}(\textbf{Comments:} #1)}}
\newcommand{\ada}[1]{#1}
\newcommand{\dla}[1]{}
\newcommand{\ria}[1]{}
\newcommand{\xo}[1]{\input{Content/#1}}
\newcommand{\oo}[1]{\input{Content/#1}}
\newcommand{\ox}[1]{}
\newtheorem{thm}{Theorem}[section]
\newtheorem{pof}{Proof}[section]
\newtheorem{defi}[thm]{Definition}
\newtheorem{prop}[thm]{Proposition}
\newtheorem{cor}[thm]{Corollary}
\newtheorem{summ}{Summary}
\newtheorem{rem}[thm]{Remark}
\newtheorem{ex}[thm]{Problem}
\newtheorem{lem}[thm]{Lemma}
\newtheorem{exa}[thm]{Example}
\newtheorem{proj}{Project}
\newtheorem{conj}{Conjecture}
\numberwithin{equation}{section}
\setcounter{tocdepth}{1}    
\hypersetup{
    colorlinks=true, 
    linktoc=all,     
    linkcolor=red,  
}

\author{Qirui Li}
\title{\papertitle}
\date{}
\maketitle
\def\load#1{\input{\papername/Command/#1}}
\def\Load#1{\input{\papername/Command/#1}}
\newcommand{\sh}{\mathfrak{h}}
\newcommand{\OOOO}[1]{\mathrm{Orb}_{#1_1,#1_2}}
\newcommand{\sg}{\mathfrak{g}}
\newcommand{\res}{\mathrm{res}}
\newcommand{\ch}{\mathrm{ch}}
\newcommand{\bc}{\mathfrak{b}}
\newcommand{\gn}{\mathrm{gl}_n}
\newcommand{\ad}{\mathrm{Ad}}
\renewcommand{\O}{\mathcal{O}}
\newcommand{\Gr}{\mathrm{Gr}}
\newcommand{\hr}{\ceil{\frac r2}}
\renewcommand{\gr}{\ceil{r}}
\newcommand{\y}[2]{\frac{#1}{1-q^{-2#2}}}
\newcommand{\Y}[3]{\frac{#1}{(1-q^{-2#2})^#3}}
\newcommand{\yy}[3]{\frac{#1}{(1-q^{-#2})(1-q^{-#3})}}
\newcommand{\tamade}{\def\zaba{\alpha_1(\zeta)}}

\newcommand{\ddl}{\dd^l}
\newcommand{\ddr}{\dd^r}

\section{Double Structures}
\lb{md}
In this paper, a double structure on an object means two different actions of a quadratic extension of $F$. 
\[defi]{Let $K$ be a quadratic etale algebra, $\OO D=\End(\GF)$ and $D=\End(\GF)\otimes_{\OO F}F$ a central simple algebra over $F$. A double $K$-structure on $\GF$ is a pair of embeddings of $F$-algebras
\[equation]{
\[split]{
\map{\alpha_1}K{\D},\\
\map{\alpha_2}K{\D}.
}
}
Two double structures $(\alpha_1,\alpha_2)$, $(\alpha_1',\alpha_2')$ are called isogenous if there is an inner automorphism $\map{c_g}\D\D$ induced by an element $g\in D$ such that the following diagrams commute
$$
\xymatrix{
&K\ar[ld]_{\alpha_1}\ar[rd]^{\alpha_1'}&\\
\D\ar[rr]^{c_g}&&\D\\
}
\qquad
\xymatrix{
&K\ar[ld]_{\alpha_2}\ar[rd]^{\alpha_2'}&\\
\D\ar[rr]^{c_g}&&\D\\
}
$$
where  $c_g(x)=g^{-1}xg$. 
The double structure is called integral if furthermore $\alpha_i(\OO K)\subset\OO D$ for $i=1,2$. Two integral double structures are called isomorphic if one can take $g\in\DDDD$ for $c_g$ in the above diagram.
}
We call this a double $K$-structure on $\GF$ because this pair gives $\GF$ two $K$-actions through quasi-isogenies. From now on, we will fix an element $\zeta\in K$ with 
$$
\zeta\notin F.
$$
It is clear that the image $\alpha_i(K)\subset D$ is determined by $\alpha_i(\zeta)$. In the rest of this paper, we denote the image of $\baso_i$ by $K_i\subset\D$. Moreover, let
$$
\dEi\subset\D,\;
\dEEi\subset\DD,\;
\dEEEi\subset\DDD,\;
\dEEEEi\subset\DDDD
$$
be centralizers of $\baso_i(\zeta)$
respectively for $i=1,2$.

The goal for this section is to establish a general theory, where we allow $\GF$ to be an arbitrary $\pi$-divisible group and we do not impose further conditions on $D$ and $K$. When we apply our general construction in this Section to Sections \S\ref{redu},\ref{indu},\ref{orba}, we all specialize to the case where $K/F$ is an unramified extension, $\GF=\left(F/\OO F\right)^{2h}$ and 
$$\D=\g\qquad\DD=\gg\qquad\DDD=\ggg\qquad\DDDD=\gggg;$$
$$\dEi=\Ei\qquad\dEEi=\EEi\qquad\dEEEi=\EEEi\qquad\dEEEEi=\EEEEi$$ 
for $i=1,2$.

Two double $K$-structures $(\biso_1,\biso_2)$ and $(\biso_1',\biso_2')$ on $\GF$ are called {\it isogenous} if there is a self quasi-isogeny $\map{\phi}\GF\GF$ carrying one double structure to another, in other words $$\phi\circ\biso_i(x)=\biso_i'(x)\circ\phi \quad \text{for any }x\in K \text{ and }i=1,2.$$ This implies $(\biso_1,\biso_2)$ and $(\biso_1',\biso_2')$ are conjugate in $D$. 

In this section, we will construct the homogeneous space $\SS$ in \eqref{huimie}, then attach an invariant polynomial to every point on $\SS$ and to every isogeny class of double structures.
 
\[defi]{\lb{ie}Let $(\alpha_1,\alpha_2)$ be a double structure, its \textit{interior angle bisector} is defined by
$$
i_{\alpha_1,\alpha_2}:=\frac{\alpha_2(\zeta)-\alpha_1(\zeta^\sigma)}{\zeta-\zeta^\sigma}.
$$
and its \textit{exterior angle bisector} is defined by
$$
e_{\alpha_1,\alpha_2}:=\frac{\alpha_2(\zeta)-\alpha_1(\zeta)}{\zeta-\zeta^\sigma}.
$$
}
\[rem]{The above definition does not depend on the choice of $\zeta\in K$. The names interior and exterior angle bisector come from the case of double structures of the Hamilton quaternion algebra where $\mathfrak{S}_2(\RR)$ is isomorphic to a sphere $S^2$. Then $\alpha_1$ and $\alpha_2$ correspond to two points on it. Joining those two points with the center of the sphere, we get an angle, and $e_{\alpha_1,\alpha_2}$ and $i_{\alpha_1,\alpha_2}$ are exactly located at the exterior and interior angle bisector respectively. Their conjugate-action are given by reflections by those bisectors.
}

\[prop]{The exterior and interior bisectors satisfy the following identities. For any $x\in K$,
$$
i_{\alpha_1,\alpha_2}\circ\alpha_1(x) = \alpha_2(x)\circ i_{\alpha_1,\alpha_2},\qquad i_{\alpha_1,\alpha_2}\circ\alpha_2(x) = \alpha_1(x)\circ i_{\alpha_1,\alpha_2},
$$
$$
e_{\alpha_1,\alpha_2}\circ\alpha_1(x) = \alpha_2(x^\sigma)\circ e_{\alpha_1,\alpha_2},\qquad e_{\alpha_1,\alpha_2}\circ\alpha_2(x) = \alpha_1(x^\sigma)\circ e_{\alpha_1,\alpha_2}
$$
and
\[equation]{\lb{chaojimiao}
i_{\alpha_1,\alpha_2}\circ e_{\alpha_1,\alpha_2} = - e_{\alpha_1,\alpha_2}\circ i_{\alpha_1,\alpha_2},\qquad (i_{\alpha_1,\alpha_2}\pm e_{\alpha_1,\alpha_2})^2=1.
}
Furthermore, we have the Pythagorean Theorem
\[equation]{\lb{lemiao}
i_{\alpha_1,\alpha_2}^2+e_{\alpha_1,\alpha_2}^2=1.
}
}
\[proof]{We have
$$(\alpha_1(x)-\alpha_2(x^{\sigma}))\alpha_2(x)=\alpha_1(x)\alpha_2(x)- \alpha_2(x^{\sigma}x).$$
Note $x^\sigma x\in F$ so $\alpha_2(x^\sigma x)=\alpha_1(x^\sigma x)$, this implies that the above equation turns into
$$
\alpha_1(x)\alpha_2(x)- \alpha_1(x^{\sigma}x)=\alpha_1(x)(\alpha_2(x)-\alpha_1(x^{\sigma}))
$$
Now since $x^\sigma+x\in F$ so $\alpha_2(x+x^\sigma)=\alpha_1(x+x^\sigma)$, we obtain
$$
\alpha_1(x)(\alpha_1(x)-\alpha_2(x^{\sigma})).
$$
Since $i_{\alpha_1,\alpha_2}=\frac{\alpha_1(x)-\alpha_2(x^{\sigma})}{x-x^\sigma}$, we proved $i_{\alpha_1,\alpha_2}\circ\alpha_2(x) = \alpha_1(x)\circ i_{\alpha_1,\alpha_2}$. The proof of the other identities are similar. The equation \eqref{chaojimiao} are followed by direct calculation and it is easy to verify $i_{\alpha_1,\alpha_2}\circ e_{\alpha_1,\alpha_2} = - e_{\alpha_1,\alpha_2}\circ i_{\alpha_1,\alpha_2}$. To prove $(i_{\alpha_1,\alpha_2} - e_{\alpha_1,\alpha_2})^2=1$, simply note that
$$
\left(i_{\alpha_1,\alpha_2} - e_{\alpha_1,\alpha_2}\right)^2 = \left(-\frac{\alpha_2(\zeta-\zeta^\sigma)}{\zeta-\zeta^\sigma}\right)^2 = 1.
$$
By expanding the left expression we also have
$$
1=\left(i_{\alpha_1,\alpha_2} - e_{\alpha_1,\alpha_2}\right)^2=i_{\alpha_1,\alpha_2}^2 + e_{\alpha_1,\alpha_2}^2 - i_{\alpha_1,\alpha_2}\circ e_{\alpha_1,\alpha_2} - e_{\alpha_1,\alpha_2}\circ i_{\alpha_1,\alpha_2}=i_{\alpha_1,\alpha_2}^2 + e_{\alpha_1,\alpha_2}^2.
$$
So we also have $1=i_{\alpha_1,\alpha_2}^2 + e_{\alpha_1,\alpha_2}^2=\left(i_{\alpha_1,\alpha_2} + e_{\alpha_1,\alpha_2}\right)^2$. We proved this proposition.
}

\[defi]{
The normalized centralizer of a double structure $(\alpha_1,\alpha_2)$ is defined by 
$$
i_{\alpha_1,\alpha_2}^2 = \frac{(\alpha_1(\zeta)-\alpha_2(\zeta^\sigma))^2}{(\zeta-\zeta^\sigma)^2}\in \dEo\cap\dEt.
$$
}
\[proof]{
\tamade
This element should be in $\D$, we prove it is in the subset $\dEo\cap\dEt
$. Since the expression is symmetric for $\alpha_1$ and $\alpha_2$, we only need to show this element commutes with $\zaba$. By the property of the interior angle bisector reflector,
$$
i_{\alpha_1,\alpha_2}\circ i_{\alpha_1,\alpha_2}\circ \zaba = i_{\alpha_1,\alpha_2}\circ \ziba\circ i_{\alpha_1,\alpha_2} = \zaba\circ i_{\alpha_1,\alpha_2}\circ i_{\alpha_1,\alpha_2}
$$
We have proved $i_{\alpha_1,\alpha_2}^2\in\dEo\cap\dEt$. }

\[defi]{\label{2.4}\comments{2.4}
The invariant polynomial of a double structure $\ds\alpha$ is the characteristic polynomial of its normalized centralizer $i_{\alpha_1,\alpha_2}^2$, as an element of the central simple algebra $\dEo$ over $K$.}

\subsection{Parameter space of double structures}
In this subsection, we define the homogeneous space $\SS$ and study the isogeny class of double structures. Given any two double structures $(\baso_1,\baso_2)$ and $(\baso_1',\baso_2')$, we may find an element $\varphi\in \DD$ such that $\baso_1=\varphi\baso_1'\varphi^{-1}$ because all embeddings $\mapp K\D$ are conjugate. Without changing its isogeny class, we may replace the pair $(\baso_1',\baso_2')$ by $(\baso_1,\varphi_1\baso_2'\varphi^{-1})$. This implies that every isogeny class of double structures has a representative $(\baso_1,-)$ with the first structure given by $\baso_1$. Then without loss of generality, we can fix an embedding $\baso_1$ and vary the second embedding $\baso_2$. Since the embedding $\map{\baso_2}KD$ is uniquely determined by the value of $\baso_2(\zeta)$, the moduli space of the embeddings $\map{\baso_2}KD$ can be represented by
\[equation]{\lb{huimie}
\SS:=\{x\in \DD: x\text{ is conjugate to }\alpha_1(\zeta)\}.
}
In other words, $\SS$ is the set of $F$-points for the conjugacy class of the matrix
$$
\m{\zeta I_h}{},{}{\zeta^\sigma I_h}.
$$
over the algebraic closure.

Our fixed $\baso_1$ determines a distinguished point $$x_0:=\baso_1(\zeta)\in \SS.$$ In the rest of this paper, we will keep the notation $x_0$. We consider the action of $\DD$ on $\SS$ via the conjugation and write $$g\cdot x := gxg^{-1}$$ for any $g\in \DD$ and $x\in\SS$. The stabilizer of $\baso_1(\zeta)$ is $\dEEo$. Moreover, any two points on $\SS$ are conjugate to each other by an element of $\DD$. This implies $\SS$ as a homogeneous space can be represented by
$$
\SS\cong \DD/\dEEo.
$$
Points on $\SS$ parametrize $F$-embeddings $\map\baso KD$. Then combining with $\baso_1$, it goes through all possible isogeny class of double structures $(\baso_1,\baso)$. Furthermore, $(\baso_1,\baso)$ and $(\baso_1,\baso')$ are isogenous if and only if $\baso$ and $\baso'$ are conjugate by an element of $\dEEo$. Therefore the space of isogeny class of double structures can be described by 
$$
\dEEo\backslash \SS\cong \dEEo\backslash \DD/\dEEo.
$$
\[defi]{For any $x\in\SS$, by the double structure induced by $x$, we mean a double structure $(\baso_1,\baso_2)$ with $\baso_2(\zeta)=x$. We also abbreviate $i_{\alpha_1,\alpha_2}$ and $e_{\alpha_1,\alpha_2}$ as $i_x$ and $e_x$. In particular,
$$
i_x = \frac{x-x_0^\sigma}{\zeta-\zeta^\sigma},\qquad e_x = \frac{x-x_0}{\zeta-\zeta^\sigma}.
$$
}

\subsection{Polar stereographic coordinate}
Suppose $F=\RR$, $K=\CC$ and $D=\HH$ a Hamilton's quaternion algebra. Then the space $\SS$ is a 2-dimensional sphere. There is a well-known universal polar stereographic coordinate system on $\SS^\circ=\SS\setminus \{x_0^\sigma\}$ in this case. In this paper we call it the \psc for short. This section is a generalization of the \psc for $\SS$ in general settings. 

\[defi]{
Consider $\D$ as a left $K$-vector space via $\baso_1$. Let $\D_+$ and $\D_-$ be eigenspaces of right multiplying $\zaba$ of eigenvalue $\zeta$ and $-\zeta$ respectively. Suppose $g$ can be decomposed as
$$
g=g_++g_-\qquad \text{for} \quad g_+\in D_+,\quad g_-\in D_-.
$$
The element
$$
\xs:=g_+^{-1}g_-
$$
is called the $\zaba^\sigma$-polar stereographic coordinate of $x=g\cdot x_0$.
}
We need to show this definition is well-defined. In other words, we need to show it only depends on $x$.
\[prop]{\lb{polarmiao}Suppose $g=g_++g_-$, then
$$
\xs=g_+^{-1}g_- = e_x\circ i_x^{-1} = (x-x_0)(x-x_0^\sigma)^{-1}.
$$
}
\[proof]{\tamade
Remember $x=\ziba$, $x_0=\zaba$ and $g$ is choosen so that $gx = x_0g$. So we have
$$
(g_++g_-)\circ\ziba = \zaba\circ(g_++g_-)
$$
Using $\zaba\circ g_-=g_-\circ\zaca$, we have
$$
g_+(\zaba-\ziba) = g_-(\ziba - \zaca)
$$
This implies
$$
g_+\circ e_x = g_-\circ i_x.
$$
We proved this proposition.}
\[rem]{The $\zaba^\sigma$-\psc can not be defined for elements $x$ such that $x_0^\sigma-x$ is not invertible.}

\renewcommand\lgggg{\mathfrak{g}_{2h}}
\renewcommand\lgg{\mathfrak{g}_{2h}}
\renewcommand\LEE{\mathfrak{g}_{a}}
\renewcommand\LEEEE{\mathfrak{g}_{a}}
\renewcommand\LEEo{\mathfrak{g}_{h,1}}
\renewcommand\LEEEEo{\mathfrak{g}_{h,1}}
\renewcommand\LEEEEt{\mathfrak{g}_{h,2}}
\renewcommand\LEEEEi{\mathfrak{g}_{h,i}}
\renewcommand\LEEt{\mathfrak{g}_{h,2}}
\renewcommand\LEEi{\mathfrak{g}_{h,i}}

\renewcommand{\LPF}{\mathfrak{p}_{F}}
\renewcommand{\LPK}{\mathfrak{p}_{K}}
\newcommand{\LDPK}{\mathfrak{p}^\Delta_{U,K}}
\newcommand{\LTPK}{\mathfrak{p}^t_{U,K}}
\renewcommand{\LPKK}{\mathfrak{p}_{U,K}}
\renewcommand{\lpK}{\mathfrak{p}_{K}}
\renewcommand{\LPKt}{\mathfrak{p}_{K,2}}
\renewcommand{\LPKi}{\mathfrak{p}_{K,i}}
\renewcommand{\lpKt}{\mathfrak{p}_{K,2}}
\renewcommand{\lpKi}{\mathfrak{p}_{K,i}}

\renewcommand{\LUF}{\mathfrak{u}_{F}}
\newcommand{\LUK}{\mathfrak u_{K}}
\newcommand{\LDUK}{\mathfrak u^\Delta_{K}}
\newcommand{\LTUK}{\mathfrak u^t_{K}}
\renewcommand{\LUKo}{\mathfrak{u}_{K,1}}
\renewcommand{\LUKt}{\mathfrak{u}_{K,2}}
\renewcommand{\LUKi}{\mathfrak{u}_{K,i}}

\renewcommand{\luf}{\mathfrak{u}_{F}}
\renewcommand{\luko}{\mathfrak{u}_{K,1}}
\renewcommand{\lukt}{\mathfrak{u}_{K,2}}
\renewcommand{\luki}{\mathfrak{u}_{K,i}}

\newcommand\AAA[2]{A_r[#1,#2]}
\newcommand\NNN{N(r)}
\newcommand\aaa[2]{a_r[#1,#2]}
\newcommand\aaz[2]{a_0[#1,#2]}
\newcommand\AAZ[2]{A_0[#1,#2]}
\newcommand\AAAA[1]{A_r(#1,X)}
\newcommand\AAAAA[1]{A_0(#1,X)}
\newcommand\aaaaa[1]{a_0(#1,X)}
\newcommand\BBB[1]{B_0[#1]}
\newcommand\CCC[2]{C_r[#1,#2]}
\newcommand\CCCC[1]{C_r(#1,X)}
\newcommand\CCCCC[2]{C_r(#1,q^{#2}X)}
\newcommand\vf[1]{\frac{\vv_F\left(\iv{#1}(1)\right)}2}

\newcommand\xs{x_\#}
\newcommand\PSC{Polar Stereographic Coordinates }
\newcommand\ag{arithmetic-geometric }
\newcommand\an{analytic }

\newcommand{\OF}{\OO F\times\OO F}
\newcommand{\OL}{\OO L}
\newcommand{\OLX}{\OO L^\times}

\newcommand{\ooo}{\one_{\gggg}}
\newcommand{\oooo}{\one_{\satset}}

\newcommand{\ir}{I_r(R_+,R_-,s)}
\newcommand{\iod}{I_r^{\text{odd}}(R_+,R_-,s)}
\newcommand{\iodm}{I_r^{\text{odd}}(R_+,R_-,-s)}
\newcommand{\iev}{I_r^{\text{even}}(R_+,R_-,s)}
\newcommand{\ievm}{I_r^{\text{even}}(R_+,R_-,-s)}
\newcommand{\ird}{I'_r(R_+,R_-,0)}
\newcommand{\ri}{I_r(R_-,R_+,-s)}
\newcommand{\rid}{I_r'(R_-,R_+,0)}
\newcommand{\irp}{I_{r+}(R_+,R_-,s)}
\newcommand{\irpd}{I'_{r+}(R_+,R_-,0)}
\newcommand{\irm}{I_{r-}(R_+,R_-,s)}
\newcommand{\irmd}{I'_{r-}(R_+,R_-,0)}
\newcommand{\jr}{J_{r}(s)}
\newcommand{\rj}{J_{r}(-s)}


\newcommand\cont[1]{\mathcal C_{#1}}
\newcommand\inte[2]{\mathcal I_{#1}^{(#2)}}

\newcommand\cen{\mathbf{w}}
\newcommand\sem{\mathbf{z}}
\newcommand\sema{\mathbf{z}_-}
\newcommand\seme{\mathbf{z}_+}
\newcommand\lat{\Lambda}
\newcommand\latp{\Lambda_+}
\newcommand\latm{\Lambda_-}
\newcommand\latz{\Lambda_0}
\newcommand\latzp{\Lambda_{0+}}
\newcommand\latzm{\Lambda_{0-}}
\newcommand\lato{\Lambda_1}
\newcommand\latop{\Lambda_{1+}}
\newcommand\latom{\Lambda_{1-}}
\newcommand\lati{\Lambda_i}
\newcommand\latip{\Lambda_{i+}}
\newcommand\latim{\Lambda_{i-}}
\newcommand\latt{\Lambda_2}
\newcommand\lattp{\Lambda_{2+}}
\newcommand\lattm{\Lambda_{2-}}
\newcommand\lath{\Lambda_3}
\newcommand\latset{\mathcal L_{\blin_1,\blin_2}}
\newcommand\satset{\mathcal U_{\blin_1,\blin_2}}
\newcommand\sitset{S_{[\Lambda_+],[\Lambda_-]}}
\newcommand\sotset{S_{[R_+],[R_-]}}
\newcommand\sutset{S_{[R_+],[R_-]}^{\text{even}}}
\newcommand\setset{S_{\text{even}}}
\newcommand\sta[1]{R_{#1}}
\newcommand\stax[1]{R_{#1}^\times}
\newcommand\stap[1]{R_{#1+}}
\newcommand\stam[1]{R_{#1-}}

\newcommand{\h}{H^{\prime\prime\prime}}
\newcommand{\gth}{\mathrm{G_{2h}}}
\newcommand{\gfh}{\mathrm{G_{2h}^\prime}}
\newcommand{\hhh}{\mathrm{H_h^\prime}}
\newcommand{\hh}{\mathrm C(\blin_1)}
\newcommand{\inv}{\mathrm{inv}}
\newcommand{\Nrd}{\mathrm{Nrd}}
\newcommand{\nrd}{\mathrm{nrd}}
\newcommand{\sat}{\mathcal{S}}
\newcommand{\pf}{\phi_F}
\newcommand{\ppf}{\Phi_F}
\newcommand{\ppk}{\Phi_K}
\newcommand{\ssf}{\mathrm{S}_F}
\newcommand{\sssf}{\mathfrak{S}_F}
\newcommand{\ssk}{\mathrm{S}_K}
\newcommand{\sssk}{\mathfrak{S}_K}
\newcommand{\Int}{\mathrm{Int}}

\newcommand{\J}{\Ker\gpp}
\newcommand{\standard}{standard}
\newcommand{\qz}[1]{{\mathcal D^{#1}}}
\newcommand{\bz}[1]{{\mathcal D_{#1}}}
\newcommand{\N}[1]{|\mathrm{Nrd}(#1)|_F^{-1}}
\tikzstyle{every picture}+=[remember picture]
\tikzstyle{na} = [shape=rectangle,inner sep=0pt,text depth=0pt]
\newcommand{\bisoo}{\biso^{-1}}
\newcommand{\blinn}{\blin^{\vee}}
\newcommand{\bbbbiso}{\Isog_{\OO E}(\GG_{E},\GF)}
\newcommand{\bbbblin}{\Isom_E(E^{lkh},F^{kh})}

\newcommand{\bbbiso}{\Isom_{\OO F}(\GK,\GF)}
\newcommand{\bbblin}{\Isom_{\OO F}(\OO K^{h},\OO F^{2h})}
\newcommand{\bilin}{\Isom_{\OO F}(\OO K^{h},\OO F^{kh})}
\newcommand{\bbiso}{\Isog_{\OO F}(\GK,\GF)}
\newcommand{\bblin}{\Isom_F(K^{h},F^{2h})}
\newcommand{\obbiso}{\Isog_{\OO F}(\GG_{K_1},\GF)}
\newcommand{\obblin}{\Isom_F(K_1^h,F^{2h})}
\newcommand{\tbbiso}{\Isog_{\OO F}(\GG_{K_2},\GF)}
\newcommand{\tbblin}{\Isom_F(K_2^h,F^{2h})}
\newcommand{\balin}{\Isom_F(K^h,F^{kh})}
\newcommand{\bibiso}{\Isog_{\OO F}(\GF,\GF)}
\newcommand{\biblin}{\Isom_F(F^{2h},F^{2h})}
\newcommand{\bublin}{\Isom_F(F^{kh},F^{kh})}
\newcommand{\boblin}{\Isom_F(F^2,K)}
\newcommand{\eqh}{\mathrm{Equi}_h(F,K)}
\newcommand{\eqth}{\mathrm{Equi}_{2h}(F,F)}
\newcommand{\heeg}{\mathrm{Heeg}_{h}(F,K)}
\newcommand{\Ltimes}{\otimes^{\mathbb L}}
\newcommand{\ms}[1]{\mm{#1_{11}}{#1_{12}}{#1_{21}}{#1_{22}}}
\newcommand{\mv}{\mm1{}{}{\pi^v}}
\newcommand{\mf}{\mm111{-1}}
\newcommand{\es}{\ep_S}
\newcommand{\er}{\ep_R}
\newcommand{\qqm}[1]{\{#1\}}
\newcommand{\qqsigma}{\sigma}
\newcommand{\qq}[2]{\qqm{#1}+\qqm{#2}\qqsigma}
\newcommand{\qspace}[1]{V_{#1}}
\newcommand{\U}[1]{U_{#1}}
\newcommand{\PD}[3]{\prod_{\vv\in{#3}}(#1-#2(\vv))}
\newcommand{\bs}{\overline{s}}

\newcommand{\biso}{\varphi}
\newcommand{\blin}{\tau}
\newcommand{\infint}[4]{c_{#1,#2}({#3}{#4})}
\newcommand{\CA}[2]{Q_{{#1},{#2}}}
\newcommand{\GGGG}[1]{\mathfrak{G}_{#1}}
\newcommand{\AC}[2]{({#1},{#2})}
\newcommand{\FO}[2]{R_{#1,#2}}
\newcommand{\GGT}[1]{\mathcal{G}[\pi^{#1}]}
\newcommand{\GGS}[1]{\mathcal{G}_0^{#1}}
\newcommand{\GGH}[3]{\GG_0^{#3}[\AC{#1}{#2}]}
\newcommand{\GGSS}[2]{\GG_0[\pi^{#1}]^{#2}}
\newcommand{\DEF}[1]{{\mathcal{M}_{#1}}}
\newcommand{\DEf}[1]{{\mathcal{N}_{#1}}}
\newcommand{\OD}[2]{D_{#1}}



\newcommand{\inter}[2]{\chi({#1}\otimes^{\mathbb{L}}{#2})}
\newcommand{\q}[3]{\newcommand{#1}{#2}
}
\newcommand{\sa}{^{\prime\prime}}
\newcommand{\Deltai}[1]{\Delta_{#1}}
\newcommand{\Deltao}{\Delta_1}
\newcommand{\Deltat}{\Delta_2}
\q{\Nrdd}{\Nrd_{D}}{Nrdd}
\q{\Nrdg}{\Nrd_{\mathbb{G}}}{Nrdg}
\q{\GG}{\mathcal{G}}{G}
\newcommand{\equi}[3]{\mathrm{Equi}_{#3}(#1/#2)}
\newcommand{\Fn}{\mathcal F_n}
\newcommand{\Fm}{\mathcal F_m}
\newcommand{\Fmv}{\mathcal{H}_{m+v}}
\newcommand{\Fnvv}{\mathcal{H}_{n+v}}
\newcommand{\pmn}{\pp{m-n}}
\newcommand{\intt}[2]{\chi\left({#1}\litimes{#2}\right)}
\newcommand{\inttt}[2]{\chi\left({#1}\litimes{#2}\right)}
\newcommand{\inttn}[2]{\intt{#1}{#2}_n}
\newcommand{\intti}[2]{\intt{#1}{#2}_\infty}
\newcommand{\jiao}[2]{\mathrm{Int}(#1,#2)}
\newcommand{\ji}[1]{\mathrm{Int}(#1)}
\newcommand{\inta}[1]{\ca{\biso_1}{\blin_1}{#1}\bullet\ca{\biso_2}{\blin_2}{#1}}
\newcommand{\hinta}[1]{\ca{\biso_1}{\blin_1}{#1}\bullet\pp{v}_*\pro^*\ca{\biso_2}{\blin_2}{#1}}
\newcommand{\pcs}[1]{\pi^{#1}\OO F^\times\oplus\mu\OO F^\times}
\newcommand{\litimes}{\otimes^{\mathbb L}}
\q{\GN}{\mathcal F_\infty^{\gpp}}{GN}
\q{\GO}{\delta_{1,\infty}}{GO}
\q{\GT}{\delta_{2,\infty}}{GT}
\q{\GI}{\delta_\infty}{GI}
\q{\ev}{\N{x\ACj^{-1}y-\pi^v\ACj^{-1}}}{ev}
\q{\pc}{\pi\OO F\oplus\mu\OO F^\times}{pc}
\q{\SPF}{\Spf\CFF q}{SPF}
\q{\fiber}{\GGH{\pi^{v}\gamma^{-1}}gh}{fiber}
\q{\GF}{\GG_F}{GF}
\q{\G}{\GG_F^{2h}}{G}
\q{\GH}{\GG_F^h}{GH}
\q{\GHO}{\GG_{F,1}^h}{GHO}
\q{\GHT}{\GG_{F,2}^h}{GHT}
\q{\GFK}{\GG_F^{kh}}{GFK}
\q{\GK}{\GG_K}{GK}
\q{\GKB}{\GG_{\overline{K}}}{GKB}
\q{\GKH}{\GG_K^{h}}{GKH}
\q{\superGKH}{\GG_K^{h}\times\overline{\GG_K^h}}{superGKH}
\q{\HG}{\mathrm{Hom}_{\OO F}(\underline{\OO F^{kh\vee}},\GF)}{HG}
\q{\HGK}{\mathrm{Hom}_{\OO K}(\underline{\OO K^{h\vee}},\GK)}{HGK}
\q{\TG}{\GF\otimes_{\OO F}\OO F^{kh}}{TG}
\q{\TGK}{\GK\otimes_{\OO K}\OO K^{h}}{TGK}
\q{\bp}{(\biso,\blin)}{bp}
\q{\bpp}{(\pi^{m}\biso,\blin)}{bpp}
\q{\bpn}{(\biso_0,\blin_0)}{bpn}
\newcommand{\gppp}[2]{(#1,#2)_{\infty}}
\newcommand{\gppi}[1]{\gppp{\biso_{#1}}{\blin_{#1}}}
\newcommand{\gppj}[2]{\gppp{\pi^{#2}\biso_{#1}}{\blin_{#1}}}
\q{\gp}{\gppp{\biso}{\blin}}{gp}
\q{\gpo}{\gppp{\biso_1}{\blin_1}}{gpo}
\q{\gpt}{\gppp{\biso_2}{\blin_2}}{gpt}
\q{\gpp}{{\gppj{}m}}{gpp}
\q{\gppo}{{\gppj1m}}{gppo}
\q{\gppt}{{\gppj2m}}{gppt}
\newcommand{\gc}{\delta}
\newcommand{\ca}[3]{\delta[#1,#2]_{#3}}
\newcommand{\cycy}[1]{\ca{\biso}{\blin}{#1}}
\newcommand{\cycw}[2]{\ca{\biso}{\pi^{#2}\blin}{#1}}
\newcommand{\cycn}[2]{\ca{\biso_{#1}}{\blin_{#1}}{#2}}
\newcommand{\cycwn}[3]{\ca{\biso_{#1}}{\pi^{#2}\blin_{#1}}{#3}}
\q{\cyio}{\ca{\biso_1}{\blin_1}{\infty}}{cyio}
\q{\cyit}{\ca{\biso_2}{\blin_2}{\infty}}{cyit}
\q{\cyc}{\cycy{n}}{cyc}
\q{\cyco}{\cycn1n}{cyco}
\q{\cyct}{\cycn2n}{cyct}
\q{\cywo}{\ca{\biso_1}{\pi^w\blin_1}{\infty}}{cywo}
\q{\cywt}{\ca{\biso_2}{\pi^w\blin_2}{\infty}}{cywt}
\q{\cy}{\cycy{\infty}}{cy}

\q{\ACj}{\gamma}{ACj}
\q{\ACx}{x}{ACx}
\q{\GE}{\GG_\eta}{GE}
\q{\GjE}{\GG_{\ACj_c\eta}}{GjE}
\q{\GEN}{{\GG_0^\prime}}{GEN}
\q{\QP}{Quotient Projection}{QP}
\q{\RP}{Restriction Projection}{RP}
\q{\LatF}{\OO F^{kh}}{LatF}
\q{\LatK}{\OO K^h}{LatK}
\q{\ACg}{g}{ACg}
\q{\ACpair}{\AC{\ACj}{\ACg}}{ACpair}
\newcommand{\pp}[1]{(\pi^{#1})}
\newcommand{\pg}[1]{(\pi^{#1})}
\q{\pro}{({\pi^{m}\ACj,\ACg_0})}{pro}
\q{\proi}{(\ACj,\ACg)_{\infty}}{proi}
\q{\proh}{({\ACj,\pi^vh\ACg})}{proh}
\q{\pron}{(\ACj_0,\pi^v\ACg_0)}{pron}
\q{\proc}{(\ACg\ACg_0,\ACj_0\ACj)}{proc}
\q{\prng}{(\ACj_0,\ACg)}{prng}
\q{\ACpairnot}{\AC{\gamma_0}{g_0}}{ACpairnot}
\q{\EndMF}{\OO {{D_F}}^{\times}\times \GL_h(\OO F)}{EndMF}
\q{\EndD}{\OO D}{EndD}
\q{\EndM}{\MM_{kh}(\OO F)}{EndM}
\q{\GLMF}{\GL_h(F)\cap\EndM}{GLMF}
\q{\dm}{\DEFf m{}h}{dm}
\q{\dwm}{\DEFf {m+w}{\eta_0}h}{dwm}
\q{\Dm}{\DEFF mF{kh}}{Dm}
\q{\djm}{\DEFf m{\ACj_c\eta}h}{djm}
\q{\dqm}{{\DEFf m{}h}^{(w)}}{dqm}
\newcommand{\Qua}[2]{[#1,#2]}
\newcommand{\infQua}[2]{(#1,#2)_{\infty}}
\newcommand{\Quaa}[1]{[\biso_{#1},\blin_{#1}]}
\newcommand{\Quasii}[1]{\Quasi_{#1}}
\q{\Quasi}{{(\pi^{m}\biso,\blin)}}{Quasi}
\q{\Quasio}{[\biso_1,\pi^{w_1}\blin_1]}{Quasio}
\q{\Quasit}{[\biso_2,\pi^{w_2}\blin_2]}{Quasio}
\q{\QQuasi}{\overline{\Quasi}}{QQuasi}
\q{\QQQuasi}{\overline{\Quasi}}{QQQuasi}
\q{\Quasin}{{[\biso_0,\blin_0]}}{Quasin}
\q{\superQuasi}{\infQua{\eta}{\pi^w(\tau\times\overline{\tau})}}{superquasi}

\newcommand{\TIMES}[1]{\!\!\!\!\underset{#1}{\times}\!\!\!\!}

\q{\HBN}{\pgDEFf{m+v+u+w}\eta hw}{HBN}
\q{\FBN}{haha}{FBN}
\q{\HSN}{\pgDEFf{m+u+w}\eta hw}{HSN}
\q{\FSN}{\GEN^h[\pi^{m+w}]}{FSN}
\q{\HBM}{\DEFF{m+u+v}F{kh}}{HBM}
\q{\FBM}{\G[{\pi^{m+v}\ACj^\prime}^{-1}\ACj^{-1}\ACg^\prime\ACg]}{FBM}
\q{\HSM}{\DEFF{m+u}F{2h}}{HSM}
\q{\FSM}{\G[{\pi^m\ACj^\prime}^{-1}\ACg^\prime]}{FSM}
\q{\RFSM}{\G[\pi^m\ACj\ACg^{-1}]}{RFSM}
\q{\BM}{\DEFF uF{2h}}{BM}
\q{\BN}{\DEFf{u}{\eta_0} h}{BN}

\q{\Mat}{M_{\eta,\tau}}{Mat}
\q{\Disc}{\Delta_{\eta,\tau}}{Disc}
\q{\Discn}{\Delta_{\eta_0,\tau_0}}{Discn}
\q{\Up}{\Gamma^{-1}}{Up}
\q{\Upn}{\Gamma^`}{Upn}
\q{\embed}{(\eta,\tau)}{embed}
\q{\MG}{\mathfrak{m}_{\G}}{MG}
\q{\MGEN}{\mathfrak{m}_{\GEN}}{MGEN}
\q{\RG}{\CFF q[[X_1,\cdots,X_{2h}]]}{RG}
\q{\RGEN}{\CFF q[[T_1,\cdots,T_{h}]]}{RGEN}

\q{\INTER}{\mathcal{F}_{m+v}\otimes{\pron^*\mathcal{F}_{m}}}{INTER}
\q{\SINTER}{s^*\mathcal{F}_{m+v}\otimes s^*{\pron^*\mathcal{F}_{m}}}{SINTER}
\q{\INTERNUM}{\inter{\mathcal{F}_{m+v}}{\pron^*\mathcal{F}_{m}}}{INTERNUM}
\q{\Tor}{\mathcal{T}or(\mathcal{F}_{m+v},\pron^*\mathcal{F}_{m})}{Tor}
\q{\INFINT}{\infint{\gamma_0}{\tau_0}{\pi^v\gamma_0^{-1}}{g_0}}{INFINT}
\q{\HM}{\mathcal{H}_m}{HM}
\q{\FM}{\mathcal{F}_{m+u}}{FM}
\q{\HV}{\mathcal{H}_{m+v}}{HV}
\q{\FV}{\mathcal{F}_{m+v+u}}{FV}

\newcommand{\AF}{{A_F}}
\newcommand{\AK}{{A_K}}
\newcommand{\AD}{{A_D}}
\newcommand{\GD}{{G_D}}
\newcommand{\GGG}{{G_?}}
\newcommand{\HF}{{H_F}}
\newcommand{\HK}{{H_K}}
\newcommand{\HD}{{H_D}}
\newcommand{\HHH}{{H_?}}
\newcommand{\TF}{{T_F}}
\newcommand{\TK}{{T_K}}
\newcommand{\TD}{{T_D}}
\newcommand{\OOF}{{\mathrm{Orb}_F}}
\newcommand{\OOK}{{\mathrm{Orb}_K}}
\newcommand{\OOL}{{\mathrm{Orb}_L}}
\newcommand{\ooF}{{\mathrm{orb}_F}}
\newcommand{\ooK}{{\mathrm{orb}_K}}
\newcommand{\OOD}{{\mathrm{Orb}_D}}
\newcommand{\OOO}{{\mathrm{Orb}_?}}
\newcommand{\OBF}{{\mathbb{O}^{reg}_F}}
\newcommand{\OBK}{{\mathbb{O}^{reg}_K}}
\newcommand{\OBB}{{\mathbb{O}^{reg}_?}}
\newcommand{\Res}{\mathrm{Res}}
\newcommand{\orb}{\mathrm{orb}}
\newcommand{\Orb}{\mathrm{Orb}}

\newcommand{\nm}[1]{#1\overline{#1}}

\newcommand{\sha}{\mathcal{S}}
\newcommand{\bi}{\mathcal{B}}
\newcommand{\cao}{\mathcal{C}}
\newcommand{\tao}{\mathcal{T}}
\newcommand{\rug}{\mathcal{R}}
\newcommand{\deri}[2]{\left.\frac{\dd}{\dd #1}\right|_{#1=#2}}

\newcommand\eso{S_{R_+,R_-}^{\mathrm{odd}}}
\newcommand\ese{S_{R_+,R_-}^{\mathrm{even}}}
\newcommand\cso{C^{\mathrm{odd}}}
\newcommand\csoo[1]{C_{#1}^{\mathrm{odd}}}
\newcommand\cse{C^{\mathrm{even}}}
\newcommand\csee[1]{C_{#1}^{\mathrm{even}}}
\newcommand\esp[3]{C_{#1,#2}^{(#3)}}
\newcommand\esm[3]{C_{#1,#2}^{(-#3)}}

\def\load#1{\input{\papername/#1}}
\def\Load#1{\input{\papername/#1}}
\begin{abstract}
Let $K/F$ be an unramified quadratic extension of a non-Archimedean local field. In a previous work \cite{qirui2017love}, we proved a formula for the intersection number on Lubin-Tate spaces. The main result of this article is an algorithm for computation of this formula in certain special cases. As an application, we prove the linear Arithmetic Fundamental Lemma for $\GL_4$ with the unit element in the spherical Hecke Algebra.

\end{abstract}
\tableofcontents
\def\use#1{}
\def\s#1{}
\def\Road#1{}
\section{Introduction}
\subsection{Motivation}
In this paper, we give an algorithm to compute intersection numbers of CM cycles in Lubin--Tate spaces in some special cases by following an explicit formula in \cite{qirui2017love}. Our goal is to identify these intersection numbers with the values of the first derivative of certain orbital integrals. This identity is known as the linear Arithmetic Fundamental Lemma (linear AFL)conjecture, and an application of our algorithm is to prove the conjecture for $\GL_4$. As we noted in the introduction part of \cite{qirui2017love}, the global motivation for the linear AFL arises from a generalization of the arithmetic Gan--Gross--Prasad conjectures proposed by Zhang \cite{zhang2017twisted}. 

\subsection{The linear AFL} We call the two sides of the linear AFL identity the \ag side and the \an side respectively. We briefly describe the objects appearing on the two sides. Let $F$ be a non-Archimedean local field with ring of integers $\OO F$. Let $\pi$ be a uniformizer of $\OO F$ and denote the residue field by $\FF_q\cong \OO F/\pi$. 
On the \ag side, we consider a 1-dimensional formal $\OO F$-module $\GF$ of height $2h$ over $\CFF q$. Let $K/F$ be an unramified quadratic extension. Choose two embeddings 
\[equation]{\label{shazishazi}\[split]{
\map{\biso_1}{\OO K}{\End(\GF)},\\
\map{\biso_2}{\OO K}{\End(\GF)}.
}
}
Each embedding $\biso_i$ gives rise to a special cycle $Z(\biso_i)$ on the Lubin--Tate space $\DEF{\GF}$ of $\GF$. The quantity of the \ag side is the intersection number $$\mathrm{Int}(Z(\biso_1),Z(\biso_2))$$ of these two cycles.

On the \an side, we consider two embeddings of $\OO F$-algebras
$$
\map{\blin_1}{\OO F\times \OO F}{\gl_{2h}(\OO F)},
$$
$$
\map{\blin_2}{\OO F\times \OO F}{\gl_{2h}(\OO F)}.
$$
Abbreviate the symbol $\GL_{2h}$ by $\mathbf{G}_{2h}$. Let $\ctri\tau\subset \gg$ be the centralizer of the image of $\tau_i$ for each $i=1,2$. We fix an element $g\in\gg$ such that 
\[equation]{\lb{hashi}\blin_2(x)=g^{-1}\blin_1(x)g\qquad \text{ for any }x\in\OO F\times\OO F.}
Using an isomorphism $\mathrm C(\blin_i)\cong \gga$, we can write any $x\in \mathrm C(\blin_i)$ as $x=(x_1,x_2)$ for $x_1,x_2\in\mathbf{G}_h(F)$. Moreover, define
$$
|x|:=\left|{\det(x_1^{-1}x_2)}\right|_F\qquad \eta_{K/F}(x):=\eta_{K/F}(\det(x_1x_2))
$$
where $\eta_{K/F}$ is the quadratic character of $K/F$.  
Let $\map f{\gg}{\RR}$ be an arbitrary smooth test function with compact support. We associate $\blin_1,\blin_2$ with the following relative orbital integral defined by
\begin{equation}\lb{xiaotina}
	\OOOO\blin(f,s):=\int_{{\mathrm C(\blin_1)\cap \mathrm C(\blin_2)}\backslash{\hh\times\hh}}f(u_1^{-1}g u_2)\eta_{E/L}(u_2)\left|u_1u_2\right|^s\dd u_1\dd u_2
\end{equation}
where we view ${\mathrm C(\blin_1)\cap \mathrm C(\blin_2)}$ as a subgroup of $\hh\times\hh$ via the diagonal embedding.  
The linear AFL conjecture states that
\[equation]{\label{wuliao}
\pm(2\ln q)^{-1}\left.\frac{\dd}{\dd s}\right|_{s=0}\OOOO\blin\left(\one_{\gggg},s\right)=\mathrm{Int}(Z(\biso_1),Z(\biso_2))
}
is a valid equation when $(\biso_1,\biso_2)$ matches with $(\blin_1,\blin_2)$ and the sign $\pm$ is chosen so that the quantity is positive. By definition, $(\biso_1,\biso_2)$ {\it matches} with $(\blin_1,\blin_2)$, if there is an isomorphism $\mapp{\End(\GF)\otimes_{\OO F}C}{\g\otimes_FC}$ between $C$-algebras for the algebraic closure $C$ of $F$ that makes the following diagram commute for $i=1,2$.
$$
\xymatrix{
	(F\times F)\otimes_F C\ar[r]^{\cong}\ar[d]_{\blin_i\otimes\id_C}&K\otimes_F C\ar[d]^{{\biso_i\otimes\id_C}}\\
	\g\otimes_FC\ar[r]&\End(\GF)\otimes_{\OO F}C
}
$$
The identity \eqref{wuliao} is conjectured to hold in more general settings if we replace $f$ in the \an side by an arbitrary spherical Hecke function on $\gg$ and correspondingly replace $Z(\biso_2)$ in the \ag side by $h_{f*}Z(\biso_2)$ via the Hecke correspondence $h_{f}:\DEF\GF\leftarrow\Gamma\rightarrow \DEF\GF$ defined by $f$. In this article, we only study the case when the test function is $\one_{\gggg}$.

\subsection{Classification of Double Structures}

In the linear AFL, our parameter is a pair of embeddings from a quadratic etale algebra $K$ to a central simple algebra $D=\End(\GF)\otimes_{\OO F}F$ over $F$. We call it a double $K$-structure on $\GF$ because this pair gives $\GF$ two $K$-actions through self-quasi-isogenies. We call $(\biso_1,\biso_2)$ an {\it integral} double $K$-structure if $\biso_i(\OO K)\subset\OO D$ for $i=1,2$. Two double $K$-structures $(\biso_1,\biso_2)$ and $(\biso_1',\biso_2')$ on $\GF$ are called {\it isogenous} if there is a self quasi-isogeny $\map{\ACj}\GF\GF$ carrying one double structure to another $$\ACj\circ\biso_i(x)=\biso_i'(x)\circ\ACj \quad \text{for any }x\in K \text{ and }i=1,2.$$ In other words, this means that the pairs $(\biso_1,\biso_2)$ and $(\biso_1',\biso_2')$ are conjugate in $D$. 

Note that both sides of the linear AFL depend only on the isogeny class of the corresponding double structures. For any quadratic etale algebra $K$, let $\zeta\in K$ such that $\zeta\notin F$. Let $\zeta^\sigma$ be its conjugate. Using the element $\zeta$, we define the invariant polynomial for a double structure $\map{\biso_i}KD(i=1,2)$ to be the characteristic polynomial of
$$
\frac{\left({\biso_1(\zeta)-\biso_2(\zeta^\sigma)}\right)^2}{(\zeta-\zeta^\sigma)^2}\in\mathrm{C}(\biso_1)\cap\mathrm{C}(\biso_2)
$$
as an element of $\mathrm{C}(\biso_1)$, which is a central simple algebra over $K$. Clearly this element does not depend on the choice of $\zeta$.

\subsection{Main results of the paper}

Our formula in \cite{qirui2017love} simplifies the \ag side and reduces the conjectural linear AFL to the following identity.
\begin{conj}Let $\map f{\gg}{\RR}$ be a spherical Hecke function, and $\ACj=(\blin_1,\blin_2)$ is a double $F\times F$-structure on $\gg$. Suppose that $(\blin_1,\blin_2)$ matches to a double $K$-structure on a division algebra $D$ of invariant $\frac1{2h}$. Let $\iv\ACj$ be the invariant polynomial of $(\blin_1,\blin_2)$. Let $\map{\alpha}{K}{\gg}$ be a map of $F$-algebras, $\iv g$ the invariant polynomial of the double $K$-structure $(\baso,g^{-1}\circ\baso\circ g)$ on $\gg$ for any $g\in\gg$. Then we have
\begin{equation}\lb{butanlianai}
\pm(2\ln q)^{-1}\left.\frac{\dd}{\dd s}\right|_{s=0}\OOOO\blin\left(f,s\right) = \frac{\ep_{F,2h}}{\ep_{K,h}^2}\int_{\gg}f(g)\Big|\Res(P_\ACj,P_\ACg)\Big|_F^{-1}\dd g
\end{equation}
where constants $\ep_{F,2h}$ and  $\ep_{K,h}$ are densities of invertible matrices in $\ggg$ and $\EEE$. The symbol $\mathrm{Res}$ represents the resultant of two polynomials.  
\end{conj}	

Our main result is a computational method for calculating the \ag side for $f=\one_{\gggg}$. As an application, we proved the identity \eqref{butanlianai} for $h=2$.
\[thm]{The equation \eqref{butanlianai} holds for $h=2$, $f(g)=\one_{\gggg}(g)$.
	}
For higher $h$, both sides of \eqref{butanlianai} are computable when we impose the following condition:
\[itemize]{
\item (*)The valuation $\vv_F(P_\ACj(1))$ is odd and coprime to h.
	}In this paper, our algorithm allows us to compute all intersection numbers for higher $h$ in the case of (*). There is also an inductive formula for orbital integrals but it seems too complicated to be practically useful. In particular, we have not succeeded identifying the inductive formulas for the two sides, except for some lower rank cases.

Now we give more details of our computational methods. The computation for the \ag side is described as follows. We see that the integrand in \eqref{butanlianai} is invariant under the action of $\EE$. Then we only need to compute the intersection number by integrating certain function over the homogeneous space $\SS=\gg/\EE$. Then we divide $\SS$ into a disjoint union of subsets with two properties. Firstly, each subset is invariant under the action of $\EEEE$. Secondly, when we have condition (*), our integrand is a constant on each subset. This method gives us an inductive formula for computing the intersection number. Finally, we prove the $h=2$ case of the linear AFL by comparing the result of computation at the end of Section \S\ref{orba} and Section \S\ref{chetui}. Our method for the \an side in \S\ref{orba} is counting lattices, which is approachable when $h=2$ since there is an easy classification of $\OO F$-subalgebras of a quadratic field extension $K$ over $F$.


This paper starts with Section \S\ref{md} to discuss double structures, which are parameters in the linear AFL identity. The constructions and lemmas in \S\ref{md} will be used repeatedly in our computation in the \an side(Section \S\ref{orba}) and the \ag side (Section \S\ref{hongzha} to Section \S\ref{chetui}). Section \S\ref{indu} gives a complete list of inductive formulae to compute the \ag side with condition (*). The calculation for $h=2$ case is done in Section \S\ref{chetui}. Section \S\ref{hongzha} and \S\ref{redu} are preparations for Section \S\ref{indu}. The reader may skip those two sections if they are willing to accept the formula \eqref{miaoya}.  

\Road{Introduction}
\Road{Abstract}
\Road{Acknowledgement}
\Road{Command/script}
\Road{etc}

\subsection{Invariant Polynomials}
In \cite{qirui2017love}, we defined the invariant polynomial by a different way than Definition \ref{2.4}. In the next proposition, we show those definitions are equivalent.
\[prop]{Let $\baso_2=g^{-1}\baso_1g$ and decompose $g=g_++g_-$ according to the decomposition $D=D_+\oplus D_-$. Then the normalized centralizer of $\ds\baso$ is given by
$$
i_{\alpha_1,\alpha_2}^2=(g_++g_-)^{-1}g_+(g_+-g_-)^{-1}g_+.
$$}
\[proof]{It suffices to prove this proposition for a Zariski dense subset. Then we may assume $g_+$ is invertible. Then $\xs=g_+^{-1}g_-$ is well defined and we have 
$$
(g_++g_-)^{-1}g_+(g_+-g_-)^{-1}g_+=	(1-\xs)^{-1}(1+\xs)^{-1}.
$$
By Proposition \ref{polarmiao}, we have
$$
(1-\xs)^{-1} = i_x\circ (i_x-e_x)^{-1},\qquad (1+\xs)^{-1} = i_x\circ(i_x+e_x)^{-1}.
$$
Therefore using identities in \eqref{chaojimiao},
\[equation*]{\[split]{
(1-\xs)^{-1}(1+\xs)^{-1} = i_x\circ (i_x-e_x)^{-1}\circ i_x\circ(i_x+e_x)^{-1} = i_x\circ (i_x+e_x)^{-2}\circ i_x=i_x^2.
}}
Therefore the proposition follows.}

\section{Integration in homogeneous spaces}\lb{hongzha}
We will compute our intersection number by an integral over $\SS$. In this section, we make a preparation by introducing some general theory of integration over homogeneous spaces. This section has 3 parts. The first two parts consists of some basic definitions. The subsection \S \ref{dasha} is an outline of our main strategy of integration in this paper. The Theorem \ref{fansile} will be used to prove Theorem \ref{wodemaya} and Proposition \ref{dongwu}. 

\subsection{Invariant Measure on Homogeneous spaces}
Let $G$ be an algebraic group over a local field $F$. We assume that the Lie algebra of $G$ is a finite-dimensional vector space over $F$. Let $S$ be a $G$-homogeneous space with a fixed base point $x_0\in S$. Let $H=\Stab\;x_0\subset G$ be the stabilizer of $x_0$. Then we have a canonical isomorphism $S\cong G/H$. In general, $S$ may not have a $G$-invariant Haar-measure. For example, the projective space $\PP^1_{\RR}$ is a $\GL_2(\RR)$-homogeneous space with no $\GL_2(\RR)$-invariant Haar-measure. It is well-known that the $G$-invariant measure exists for $S=G/H$ if and only if their modular characters $\delta_G$ and $\delta_H$ satisfies $\delta_G(h)=\delta_H(h) \text{ for any }h\in H$. For the rest of the paper, we only consider the case where $H$ is a compact subgroup. Then any character from $H$ to $\RR^{\times}_{>0}$ is trivial. This fact implies that we have a $G$-invariant measure on  $S\cong G/H$.

Once we have chosen a left Haar-measure $\dd g$ on $G$ and a Haar-measure $\dd h$ on $H$, the Haar-measure $\dd s$ on $S$ is defined so that for any function $\map fG\RR$, we have
$$
\int_G f(g)\dd g=\int_S\widetilde f(s)\dd s
$$
where $\widetilde f(s)=\int_H f(g_sh)\dd h$ with $g_sx_0=s$.
\subsection{Standard Haar-measure}
For an algebraic group $G$ over $\OO F$ and a smooth $G$-homogeneous space $S$, we may choose a Haar-measure on $S$ such that the total volume of $S(\OO F)$ is given by
$$
\mathrm{Vol}(S(\OO F))=\frac{\# S(\FF_q)}{q^{\dim(S)}}.
$$

The standard Haar-measure on $G$ is defined in the same way.
\subsection{Integration by fibration}\lb{dasha} In this subsection, we introduce the main computational strategy in our paper. Theorem \ref{fansile} in this section can be used in situations where the following three conditions are satisfied.

\noindent\textbf{Condition 1 (Fibration):} There is a subgroup $C\subset G$, a $C$-invariant subset $S^\circ \subset S$, a $C$-homogeneous space $T$, and a $C$-equivariant surjective map
\[equation]{\lb{tamade}
\map p{S^\circ}T.
}

\noindent\textbf{Condition 2 (Fiber Translation):} For any $t\in T$, there exists a subgroup $P_t\subset G$ such that each fiber $p^{-1}\{t\}$ is a subset of a $P_t$-homogeneous space $R_t = P_t\cdot t\subset S$. We denote $R_t^\circ = p^{-1}\{t\}$. 

\noindent\textbf{Condition 3:} We have $\dim(R_t)+\dim(T)=\dim(S)$. In other words, we require the fiber $p^{-1}\{t\}$ and $R_t$ have the same dimension.

We give an example of a fibration with those 3 conditions. If we want to project a sphere $S^2$ to its equator along longitude lines, we can not do so because of the north and south poles. Therefore, we must choose a subset $S^\circ = S^2\setminus\{\text{poles}\}$ and apply the projection. Let $t$ be a point on the equator. Each fiber $p^{-1}\{t\}$ of the projection is a semi-circle, which is a subset of a full circle $R_t\subset S^2$ passing through two poles and $t$. The full circle is a homogeneous space for a subgroup $\mathrm{O}_2(\mathbb R)\cong P_t\subset\Aut(S^2)\cong\mathrm{O}_3(\RR)$. We see the dimension of the fiber $p^{-1}\{t\}$ and $R_t$ are the same.

When we have the above conditions, we may write the integral over $S^\circ$ into the following form
\[equation]{\lb{jianbi}
\int_{S^\circ} f(s)\dd s=\int_{T}\int_{R_t^\circ}f(r) |J_t(r)|_F\dd r\dd t.
}
Here $\dd r$, $\dd t$, $\dd s$ are standard $G$,$C$,$P$-invariant measures. The main obstacle is computing the Jacobian determinant $J_t(r)$. 

We need to introduce more notation before Theorem \ref{fansile}. For any vector space $V$, by $\bigwedge V$ we mean the highest wedge product of $V$. It is an one-dimensional vector space. Note that whenever we have an exact sequence 
$$
\xymatrix{0\ar[r]& U\ar[r]&V\ar[r]&W\ar[r]&0},
$$
we have $\bigwedge V = \bigwedge U\otimes \bigwedge W$. For any $\lambda_v\in \bigwedge V$ and $0\neq \lambda_u\in \bigwedge U$, by $\frac{\lambda_v}{\lambda_u}$ we mean the unique element $\lambda_w\in\bigwedge W$ such that $\lambda_v=\lambda_w\otimes\lambda_u$.
\[thm]{\lb{fansile}Let $r\in S$ and $t=p(r)$ be as in \eqref{tamade} and \eqref{jianbi}. Moreover let $H=\Stab\;r$ be the stabilizer of $r\in S\cong G/H$. Let $C,P_t\subset G$ be subgroups as we defined above. Let $\mathfrak g$, $\mathfrak c$, $\mathfrak p$, $\mathfrak h$ be Lie algebras of $G$,$C$,$P_t$,$H$ respectively.  Let
$$
\mathfrak u=\frac{\mathfrak g}{\mathfrak p}\qquad \mathfrak u_c=\frac{\mathfrak c}{\mathfrak c\cap\mathfrak p}\qquad \mathfrak u_h=\frac{\mathfrak h}{\mathfrak h\cap\mathfrak p}
$$
be their quotient spaces.  Let

$$\dd g\in\bigwedge \mathfrak g^\vee\qquad\dd h\in\bigwedge\mathfrak h^\vee\qquad\dd p\in\bigwedge \mathfrak p^\vee\qquad\dd q\in\bigwedge(\mathfrak p\cap\mathfrak h)^\vee\qquad \dd t\in\bigwedge \mathfrak u_c^\vee$$
be volume forms corresponding to left Haar-measures on $G$, $H$, $P$, $P\cap H$ and $C/C\cap P$ respectively. Furthermore, let
$$
\dd u = \frac{\dd g}{\dd p}\in\bigwedge\mathfrak u^\vee\qquad \dd u_h=\frac{\dd h}{\dd q}\in\bigwedge\mathfrak u_h^\vee.
$$
Suppose $\dim(S)=\dim(T)+\dim(P_t\cdot r)$. Then we have $\mathfrak u=\mathfrak u_c\oplus\mathfrak u_h$ and
$$
\dd u = J(r)\dd u_h\dd t.
$$
}
\begin{proof}Let $P_tr\subset S$ be the orbit of $r$ under the action of $P_t$. Then we have $P_tr\cong P_t/P_t\cap H$. The tangent spaces of the point $r\in S^\circ$ in $Pr$ and in $S$ are naturally isomorphic to $\frac{\mathfrak g}{\mathfrak h}$ and $\frac{\mathfrak p}{\mathfrak p\cap\mathfrak h}$ respectively. Note that the natural inclusion $P_t/H\cap P_t\subset S$ gives rise to an inclusion of tangent spaces $\frac{\mathfrak p}{\mathfrak p\cap\mathfrak h}\subset \frac{\mathfrak g}{\mathfrak h}$. Since $\map p{S^\circ}T$ is an open map, the quotient of $\frac{\mathfrak g}{\mathfrak h}$ by $\frac{\mathfrak p}{\mathfrak p\cap\mathfrak h}$ is naturally isomorphic to the tangent space of $t$ in $T$. Let $T_t$ be the tangent space of $t\in T$.  The map $\mapp {\mathfrak c}{T_t}$ factors through $\mathfrak u_c$. Since $\map p{S^\circ}T$ is surjective, the induced map $$\map{p_{t,*}}{\mathfrak u_c}{T_t}$$ is naturally surjective. Given that $\dim(\mathfrak u_c)\leq \dim(S)-\dim(P_t\cdot r)=\dim(T)=\dim(T_t)$, we know that $p_{t,*}$ is an isomorphism. Therefore we have the following exact sequence 
\[equation]{\lb{ds}
\xymatrix{
	0\ar[r]&\frac{\mathfrak p}{\mathfrak p\cap \mathfrak h}\ar[r]&\frac{\mathfrak g}{\mathfrak h}\ar[r]&\mathfrak u_c\ar[r]&0.
	}
}

To prove $\mathfrak u=\mathfrak u_c\oplus\mathfrak u_h$, we complete the above sequence into the following exact diagram.
$$
\xymatrix{
	&0\ar[d]&0\ar[d]&&\\
	0\ar[r]&\mathfrak p\cap \mathfrak h\ar[r]\ar[d]&\mathfrak h\ar[r]\ar[d]&\mathfrak u_h\ar[r]&0\\
	0\ar[r]&\mathfrak p\ar[r]\ar[d]&\mathfrak g\ar[r]\ar[d]&\mathfrak u\ar[r]&0\\
	0\ar[r]&\frac{\mathfrak p}{\mathfrak p\cap \mathfrak h}\ar[r]\ar[d]&\frac{\mathfrak g}{\mathfrak h}\ar[r]\ar[d]&\mathfrak u_c\ar[r]&0\\
	&0&0&&
	}
$$
By the Snake Lemma, we have the following exact sequence
$$
\xymatrix
{
	0\ar[r]&\mathfrak u_h\ar[r]&\mathfrak u\ar[r]&\mathfrak u_c\ar[r]&0.
	}
$$
The natural inclusion $\mathfrak c\subset \mathfrak g$ induces an inclusion $\mathfrak u_c\subset\mathfrak u$. This implies that the above exact sequence splits. Then $\mathfrak u=\mathfrak u_c\oplus\mathfrak u_h$. To compute the Jacobian determinant, we note that $J(r)$ is defined in the following equation
$$
J_t(r)\dd r\dd t=\dd s.
$$
Since we have $\dd s=\frac{\dd g}{\dd h}$, $\dd r=\frac{\dd p}{\dd q}$. We may write
$$
J_t(r)\frac{\dd p}{\dd q}\cdot\dd t=\frac{\dd g}{\dd h}\in\bigwedge\left(\frac{\mathfrak g}{\mathfrak h}\right)^\vee.
$$
By multiplying $\frac{\dd h}{\dd p}$, this equation can be written to
$$
J_t(r)\frac{\dd h}{\dd q}\cdot\dd t=\frac{\dd g}{\dd p}\in\bigwedge\mathfrak u^\vee.
$$
This implies
$$
J_t(r)\dd u_h\cdot\dd t=\dd u.
$$
This completes the proof.
\end{proof}

\section{Parabolic reduction formula}
\lb{redu}

This section provides an integration formula that will be used in Section \S\ref{indu}. Let $\SSS=\SS\cap\gggg$, and $a$ a non-negative integer. For any $t$, let $P_t$ denote the corresponding invariant polynomial. Denote
$$
\SSSSSS rha=\left\{t\in\SSS: P_t(X) \text{ has exactly } a \text{ factors }(X-\lambda) \text{ with }\eico\right\}.
$$
Furthermore, we put $\SSShr {a}r:=\SSSSSS raa$ and $\SSShrr ar:=\SSSSSS ra0$. Let $\map f{P_h(F)}{\RR}$ be a function on the set $P_h(F)$ of degree $h$ monic $F$-polynomials.

We will show that
\[equation]{\lb{miaoya}
\frac1{\ep_{K_h}}\int_{\SSSSSS rha} f(P_t)\dd t=\frac1{\ep_{K_{h-a}}}\frac1{\ep_{K_a}}\int_{\SSShr ar}\int_{\SSShrr {h-a}r}f(P_{t_1}P_{t_2})|\Res(P_{t_1},P_{t_2})|_F\dd t_1\dd t_2.
}

This formula is the key for our induction algorithm.
The whole section is a proof for this. The reader willing to accept this formula may skip this section.

This section will follow the strategy of integration-by-fibration as we mentioned in Section \S\ref{hongzha}. In subsection \S\ref{fg} we introduce our basic construction of fibration and check that it satisfies the three conditions in Section \S\ref{dasha}. Then in subsection \S\ref{xiaobiaoza}, we use our constructions in Section \S\ref{md} to compute the corresponding Jacobian determinant.

\subsection{Fibration over Grassmannian}\lb{fg}
Now we construct a fibration that satiesfies the Condition 1. Consider $V=F^{2h}$ as a free $K_1$-module of rank $h$ by the $K_1$-structure induced by $\baso_1$. Let $\Gr$ be the Grassmannian variety parametrizing $n$-dimensional $K_1$-subspaces in $V$. For any $x\in\SSSSSS rha$, we may decompose the invariant polynomial of $x$ into
\[equation]{\lb{wabi}
\iv x= \ivv x\ivvv x
}
where $\ivv x$ is the maximal factor of $\iv x$ such that all roots $\lambda$ of $\ivv x$ satisfies $\eico$. Let $U_x$ be the image of the operator $\ivv x\left(i_x^2\right)$. Then $U_x$ is the maximal invariant subspace of $i_x^2$ such that all its eigenvectors in $U_x$ have its eigenvalue $\lambda$ satisfying $\eico$.
We define a map by
\[equation]{\lb{renleimiejue}
\maps p{\SSSSSS rha}\Gr{x}{U_x}.
}
Clear $\SSSSSS rha$ is a $\EEEEo$-invariant subset. This map is $\EEEEo$-equivariant.

Next we check Condition 2, that for any $U\in\Gr$, each fiber $p^{-1}(U)$ of the map $p$ in \eqref{renleimiejue} is a subset of certain homogeneous space of a certain subgroup. In our scenario here, this subgroup is $\PU\subset\gggg$, which is the stabilizer of $U$. The corresponding homogeneous space is a subset of $\SSS$ defined by 
\[equation*]{
\[split]{
\TT &:= \PF\cdot \zaba\\
&=  \{x\in\SSS: x=g\zaba g^{-1} \text{ for some }g\in\PF\}.
}
}
To check the Condition 2, we need to check the fiber is a subset of this homogeneous space.
\[lem]{\lb{inva}For any $U\in\Gr$, we have $p^{-1}(U)\subset \TT$}.
\[proof]{Let $x\in p^{-1}(U)$, we have $i_x^2U=U$.  Since any $x\in\SS$ commutes with the element $i_x^2$, the subspace $xU$ is also an invariant subspace for $i_x^2$. Furthermore, the restriction of $i_x^2$ on $U$ and $xU$ gives the same eigenvalues. Since $U$ is the maximal invariant subspace with eigenvalue $\lambda$ of $i_x^2$ that satisfies $\eico$, then we must have
$$
xU\subset U.
$$
This implies $xU=U$. Similarly, we have $\zaba U=U$. This implies 
$$
(\zeta-\zeta^\sigma)i_x = x-\zaba^\sigma\in\PU.
$$
Now
$$
x=i_x\zaba i_x^{-1}\in \TT
$$ and we are done.}

Now we check the Condition 3. As varieties over $\OO F$, we have $\dim(\SSS)=2h^2$, $\dim(\TT)=2a^2+2(h-a)^2+2a(h-a)$, $\dim(\Gr)=2a(h-a)$. Then we have $\dim(\SSS)=\dim(\TT)+\dim(\Gr)$. 

\subsection{Computation of the Jacobian Determinant}\lb{xiaobiaoza}
From now on, we denote the fiber $p^{-1}(U)$ of $U$ by $\TTT$. Let $\dd U$(resp.$\dd t'$,$\dd t$) be the $\EEEEo$(resp. $\PK$,$\gggg$)-invariant standard Haar-measure on $\Gr$ (resp. $\TT$,$\SS$).  We can write 

\[equation]{\lb{shaleba}
\int_{\SSSSSS rha}f(P_t)\dd t = \int_{\Gr}\int_{\TTT}f(P_{t'})|J(t)|_F\dd t'\dd U.
}

The Jacobian determinant $J(t)$ 
is defined by
$$
\dd t = J(t)\dd t'\dd U.
$$
We determine the value of $J(t)$ via the following theorem.

\[thm]{\lb{wodemaya}Let $x\in \SSSSSS rha$ be an element whose invariant polynomial $\iv x$ decomposes as
$$
\iv x = \ivv x\ivvv x
$$
as in \eqref{wabi}. Then we have
$$
|J(x)|_F=|\mathrm{Res}(\ivv x,\ivvv x)|_F.
$$
}

\begin{proof}
Let $\map{\baso_2} K{F^{2h}}$ be the map with $\ziba=x$. Denote the stabilizer of $x\in\SSS$ by $\EEEEt$. Let $\lgggg$, $\LEEEEo$, $\LEEEEt$, $\LPF$, $\LPKi$ be Lie algebras for groups 
$\gggg$, $\EEEEo$, $\EEEEt$, $\PF$, $\PKi$ respectively for $i=1,2$. Denote the standard Haar-measure for them as $\dd g$,$\dd k_1$,$\dd k_2$,$\dpf$,$\dpi$ respectively. For $i=1,2$, we denote
$$
\LUKi:=\frac{\LEEEEi}\LPKi \qquad \LUF:=\frac{\lgggg}{\LPK}.
$$
Let $\dd u_1=\frac{\dko}{\dpo}$, $\dd u_2=\frac{\dkt}{\dpt}$ $\dd u_F=\frac{\df}{\dpf}$. Applying Theorem \ref{fansile}, we have
\[equation]{\lb{wandan}
\dd u_F=J(x)\dd u_1\cdot\dd u_2\in \bigwedge \left(\frac{\lgggg}\LPF\right)^\vee.
}
Note that $\luki$ are tangent spaces of Grassmannian manifolds at the point of $U_x$. Let $$V=F^{2h}.$$ Recall that
$$
	U_x\cong F^{2a}.
$$
	We have natural isomorphisms
$$
\luki\cong\Hom_{K_i}(U_x,V/U_x)\qquad \luf\cong\Hom_F(U_x,V/U_x)\quad i=1,2.
$$
	Before we calculate $J(t)$, let $\baso_i'$ and $\baso_i''$ be the induced $K_i$-structures on $U_x$ and $V/U_x$. Let $\overline{\baso_i''}$ be the Galois conjugate of $\baso_i''$. Then
$(\baso_1',\baso_2')$ and $(\baso_1'',\baso_2'')$ are double structures on $U_x$ and $V/U_x$ respectively.  
Let
$$
\hipa=\dhipa\quad\hima=\dhima
$$
be subsets of $K$-linear homomorphisms with $K$-structure induced by $(\baso_i',\baso_i'')$ and $(\baso_i',\overline{\baso_i''})$ respectively for $i=1,2$.  
For any $i=1,2$, the following sequence is exact
	\[equation]{\lb{exact}
	\xymatrix{
		0\ar[r]&\hipa\ar[r]&\hof\ar[rrrrr]^(.6){T_i:= f\mapsto f\circ\baso_i'(\zeta)-\baso_i''(\zeta)\circ f}&&&&&\hima\ar[r]&0
}
}
The last map is a surjective map because its right inverse is given by
$$
	\maps{T_i'}{\hima}{\hof}f{f\circ \baso_i(\zeta-\zeta^\sigma)^{-1}\circ(1-\sigma)}_.
$$
This implies
$$
\hof = \hipa\oplus\hima.
$$
	Let $\dd u_F$,$\dd u_i$ and $\overline{\dd u_i}$be the standard Haar-measures on $\Hom_F(U_x,V/U_x)$, $\hipa$, $\hima$ respectively. By the exact sequence \eqref{exact} and the fact that $T|_{\hima}(f)=f\circ\baso_i'(\zeta-\zeta^\sigma)$, we have
$$
	\dd u_F=\det(T_i|_{\hima})\dd u_i\overline{\dd u_i}=\mathrm{Disc}_{K/F}^{a(h-a)}\dd u_i\overline{\dd u_i}.
$$
Now we start our calculation of $J(x)$. On one hand, we consider the following exact sequence
$$
\xymatrix{
	0\ar[r]&\htpa\ar[r]\ar[d]&\hopa\oplus\htpa\ar[r]\ar[d]&\hopa\ar[r]\ar[d]&0\\
	0\ar[r]&\htpa\ar[r]&\hof\ar[r]^(.7){T_2}&\htma\ar[r]&0
	}_.
$$
Then we have
$$
\dd u_2\overline{\dd u_2}=\dd u_F=J(x)\dd u_1\dd u_2\implies J(x)=\frac{\overline{\dd u_2}}{\dd u_1}
$$
This implies $J(x)$ is the relative determinant
$$
	J(x)=\det\left(\hopa\longrightarrow\hof\longrightarrow\htma\right)=\det(T|_{\hopa}).
$$
Here the map $T_2$ in the exact sequence is given by the following map
$$
\maps {T_2}\hof\htma f{f\circ \zibap-\zibapp\circ f}_.
$$

\newcommand\JJJ{J}
For any symbol $a$ we denote $L_a$ for the left composing map $L_a:f\mapsto a\circ f$ and $R_a$ for right composing map $R_a:f\mapsto f\circ a$. Please note that for any symbol $a,b$, $L_a$ and $R_b$ always commute(because composition of maps is associative). Consider the following linear operators

\[equation]{\label{minus}
\maps{\Theta_-}\hof\hof f{f\circ(\baso'_1(\zeta)-\baso'_2(\zeta))-(\baso''_1(\zeta)-\baso''_2(\zeta))\circ f,}
}
\[equation]{\label{plus}
\maps{\Theta_+}\hof\hof f{f\circ(\baso'_1(\zeta)-\baso'_2(\zeta))+(\baso''_1(\zeta)-\baso''_2(\zeta))\circ f,}
}
\[equation]{\label{inter}
\maps{L_{i''}\circ R_{e'}}\hof\hof f{i_{\alpha_1'',\alpha_2''}\circ f\circ e_{\alpha_1',\alpha_2'},}
}
\[equation]{\label{conju}
\maps{L_{i''}\circ L_{e''}}\hof\hof f{i_{\alpha_1'',\alpha_2''}\circ e_{\alpha_1'',\alpha_2''}\circ f.}
}
Then we know $\Theta_-|_{\hopa} = T_2|_{\hopa}$ and our goal is to compute
$$
\det\left(\hopa\overset{\Theta_-}\longrightarrow\homa\right).
$$
Since $U_x\cong K^a$ we have 
$$
\hipa\cong (V/U_x)^a,\qquad \hima\cong (V/U_x)^a .
$$
Therefore
$$
\det\left(\xymatrix{\hipa\ar[rr]^{L_{i''}\circ L_{e''}}&&\hima}\right) = \det\left(\xymatrix{(V/U_x)^a\ar[rr]^{L_{i''}\circ L_{e''}}&&(V/U_x)^a}\right) = \det(i_x''\circ e_x'')^a . 
$$
So we proved 
$$\det\left(\xymatrix{\hopa\ar[rr]^{L_{i''}\circ L_{e''}}&&\homa}\right) = \det\left(\xymatrix{\htpa\ar[rr]^{L_{i''}\circ L_{e''}}&&\htma}\right).$$
Use this identity and the following commutative diagram
$$
\xymatrix{
	\hopa\ar[rr]^{L_{i''}\circ L_{e''}}\ar[d]_{L_{i''}\circ R_{e'}}&&\homa\ar[d]^{L_{i''}\circ R_{e'}}\\
	\htma\ar[rr]^{L_{i''}\circ L_{e''}}&&\htpa\\
}_,
$$
we know that
$$
\det\left(\xymatrix{\hopa\ar[rr]^{L_{i''}\circ R_{e'}}&&\htma}\right) = \det\left(\xymatrix{\homa\ar[rr]^{L_{i''}\circ R_{e'}}&&\htpa}\right).
$$
Again using the exact sequence \eqref{exact}, we have the following commutative diagram
$$
\xymatrix{
	0\ar[r]&\hopa\ar[r]\ar[d]_{L_{i''}\circ R_{e'}}&\hof\ar[r]\ar[d]_{L_{i''}\circ R_{e'}}&\homa\ar[r]\ar[d]_{L_{i''}\circ R_{e'}}&0\\
	0\ar[r]&\htpa\ar[r]&\hof\ar[r]&\htma\ar[r]&0
}_.
$$
This implies $\det\left(\xymatrix{\hof\ar[rr]^{L_{i''}\circ R_{e'}}&&\hof}\right)$ equals to 
$$
\det\left(\xymatrix{\hopa\ar[rr]^{L_{i''}\circ R_{e'}}&&\htma}\right)\cdot\det\left(\xymatrix{\homa\ar[rr]^{L_{i''}\circ R_{e'}}&&\htpa}\right) = \det\left(\xymatrix{\hopa\ar[rr]^{L_{i''}\circ R_{e'}}&&\htma}\right)^2 .
$$
Since the map $\xymatrix{\hof\ar[rr]^{L_{i''}\circ R_{e'}}&&\hof}$ does not change when changing indces from $1$ to $2$, we conclude that
$$
\det\left(\xymatrix{\hopa\ar[rr]^{L_{i''}\circ R_{e'}}&&\htma}\right) = \det\left(\xymatrix{\htma\ar[rr]^{L_{i''}\circ R_{e'}}&&\hopa}\right).
$$
This identity and the following commutative diagram
$$
\xymatrix{
	\hopa\ar[d]_{L_i''\circ R_e'}\ar[rr]^{\Theta_-}&&\htma\ar[d]^{L_i''\circ R_e'}\\
	\htma\ar[rr]^{\Theta_+}&&\hopa\\
	}
$$
imply that
$$
\det\left(\xymatrix{\hopa\ar[rr]^{\Theta_-}&&\htma}\right) = \det\left(\xymatrix{\htma\ar[rr]^{\Theta_+}&&\hopa}\right).
$$
Since
$$
	\Theta_+\circ\Theta_- = (L_{\baso_1''(\zeta)-\baso_2''(\zeta)} + R_{\baso_1'(\zeta)-\baso_2'(\zeta)})\circ(L_{\baso_1''(\zeta)-\baso_2''(\zeta)} - R_{\baso_1'(\zeta)-\baso_2'(\zeta)}),
	$$
we have
	$$
\Theta_+\circ\Theta_- = L_{\baso_1''(\zeta)-\baso_2''(\zeta)}^2 - R_{\baso_1'(\zeta)-\baso_2'(\zeta)}^2.
	$$
Since the characteristic polynomial of $\frac{(\baso_1''(x)-\baso_2''(x))^2}{(x-x^\sigma)^2}$ and $\frac{(\baso_1'(x)-\baso_2'(x))^2}{(x-x^\sigma)^2}$ are $\ivv x$ and $\ivvv x$ respectively, we have
$$
{\det}_K(\Theta_+\circ\Theta_-|_{\hopa}) = \mathrm{Res}(\ivv x,\ivvv x)\mathrm{Disc}_{K/F}^{a(h-a)} .
$$
Therefore $\left|\det\left(\xymatrix{\hopa\ar[rr]^{\Theta_-}&&\htma}\right)\right|_F$ equals
$$
\sqrt{|{\det}_K(\Theta_+\circ\Theta_-|_{\hopa})|_K} =\sqrt{|\mathrm{Res}(\ivv x,\ivvv x)|_K} = |\mathrm{Res}(\ivv x,\ivvv x)|_F|\mathrm{Disc}_{K/F}^{a(h-a)}|_F.
$$
This completes the proof.\end{proof}

By Theorem \ref{wodemaya}, we can write \eqref{shaleba} as
\[equation]{\lb{ash}\[split]{
\int_{\SSSSSS rha}f(P_t)\dd t &= \int_{\Gr}\int_{\TTT}f(P_{t})|\mathrm{Res}(\ivv t,\ivvv t)|_F\dd t\dd U\\
&=\Vol(\Gr)\int_{\TTT}f(P_{t})|\mathrm{Res}(\ivv{t},\ivvv{t})|_F\dd t.
}}
Since
$$
\Vol(\Gr)=\frac{\ep_{K_h}}{\ep_{K_m}\ep_{K_n}},
$$
we have
\[equation]{\lb{laiceshi}
\frac1{\ep_{K_h}}\int_{\SSSSSS rha}f(P_t)\dd t=\frac1{\ep_{K_m}\ep_{K_n}}\int_{\TTT}f(P_{t})|\mathrm{Res}(\ivv{t},\ivvv{t})|_F\dd t.
}
Let $\one_{\TTT}$ be the characteristic function of the subset $\TTT\subset\TT$. We can write
\[equation*]{\lb{leishi}\[split]{
\int_{\TTT}f(P_{t})|\mathrm{Res}(\ivv{t},\ivvv{t})|_F\dd t&=\int_{\TT}\OTTT(t)f(P_{t})|\mathrm{Res}(\ivv{t},\ivvv{t})|_F\dd t.
}}
Consider a fibration
$$
\mapi{\TT}{\SSSh a\times\SSSh{h-a}}{(\baso_1,\baso_2)}{\left((\baso_1',\baso_2'),(\baso_1'',\baso_2'')\right)}_.
$$
This is an $\PF$-equivariant map and both spaces are $\PF$-homogeneous spaces. Each fiber is a homogeneous space of its unipotent subgroup 
$$\UF=\{g\in\PF: g|_U=\id_U, g|_{V/U}=\id_{V/U}\}.$$
Let $\dd t_1,\dd t_2$ and $\dd u$ be standard Haar-measures on $\SSSh a$, $\SSSh{h-a}$ and $\UF$ respectively. We may write the integral
\[equation*]{\[split]{
	&\int_{\TT}\OTTT(t)f(P_{t})|\mathrm{Res}(\ivv{t},\ivvv{t})|_F\dd t\\
	=&\int_{\SSSh a}\int_{\SSSh{h-a}}\int_{\UF}\OTTT(t)f(P_{t})|\mathrm{Res}(\ivv{t},\ivvv{t})|_F\dd u\dd t_1\dd t_2\\
	=&\int_{\UF}\dd u\int_{\SSShr ar}\int_{\SSShrr {h-a}r}f(P_{t_1}P_{t_2})|\Res(P_{t_1},P_{t_2})|_F\dd t_1\dd t_2.
}}
Since $\Vol(\UF)=1$, we have

\[equation]{\lb{maodian}
\frac1{\ep_{K_h}}\int_{\SSSS rha} f(P_t)\dd t=\frac1{\ep_{K_{h-a}}}\frac1{\ep_{K_a}}\int_{\SSShr ar}\int_{\SSShrr {h-a}r}f(P_{t_1}P_{t_2})|\Res(P_{t_1},P_{t_2})|_F\dd t_1\dd t_2.
}

\section{Inductive formulae for the intersection number}
\lb{indu}

In this section, we introduce the inductive formulae for the intersection number $\intnu$, where the input double structure $\ACj=(\biso_1,\biso_2)$ satisfies the following condition.
\[itemize]{
\item (*)The valuation $\vv_F(P_\ACj(1))$ for the invariant polynomial at 1 is odd and coprime to h.
}
Then let $r=\frac{-\vv_F(P_\ACj(1))}h$, which also equals to $\vv_F(\ACj_\#^2)$ where $\ACj_\#$ is the \psc of $(\biso_1,\biso_2)$. This section is divided into 3 parts. 

In the first part, we simplify the intersection formula by defining the following integrals.
\begin{equation}\lb{aaa}
	\AAA ak=\frac1{\ek a}\int_{\SSShr ar}|\iv x(1)|_F^{-k}\dd x,
\end{equation}
\begin{equation}\lb{bbb}
	\BBB c=\frac1{\ek c}\int_{ \SSShrr {c}0}\dd x,
\end{equation}
\begin{equation}\lb{ccc}
	\CCC bc=\frac1{\ek b}\int_{\SSShrr {b}r\cap\SSShr {b}0}|\iv{x}(1)|^{-b-c}_F\dd x.
\end{equation}
We prove that under the condition (*) the number $\intnu$ can be written into the following form
\begin{equation}\lb{jiaojiao}
\intnu=\sum_{a+b+c=h}|\iv\ACj(1)|_F^{-a}\AAA a0\CCC bc\BBB c.
\end{equation}
In particular, in the situation of (*), the intersection number only depends on $r$. To simplify our notation in the rest of the paper, we will denote it by $\NNN:=\intnu$.

The second part will introduce the inductive formulae for computing $\AAA a0$,$\CCC bc$,$\BBB c$. Specifically, we will show 
\begin{equation}\lb{a3}
	\CCC nm=\AAZ n{n-m}-\sum_{i=0}^{n-1} \CCC im\AAA{n-i}{n-m} \text{ for }m\geq n ,
\end{equation}
\[equation]{\lb{a4}
\BBB a= \prod_{i=1}^a\frac{1-q^{1-2i}}{1-q^{-2i}}-\sum_{i=1}^a\AAZ i0\BBB{a-i}.}
The formulae \eqref{bbb},\eqref{ccc} and \eqref{jiaojiao} imply that the computation of $\intnu$ can be reduced to $\AAA a0$. 

In the third part of this section, we will use polar stereographic coordinates to prove the following formula. Define
\begin{equation}\lb{qusi}
	\aaz n{n-m}:=q^{-2m}\prod_{\substack{i=0\\i\neq m}}^{n-1}\frac{1}{1-q^{-2(m-i)}},
\end{equation}
and define recursively
\begin{equation}\lb{a1}
\aaa{n}{n-m}:=\aaz n{n-m}-\sum_{i=1}^{m}\CCC im\aaa{n-i}{n-m}\text{ for }0\leq m<n.
\end{equation}
Then we claim
\begin{equation}\lb{a2}
\AAA n{n-m}=\sum_{i=0}^{n-1} \aaa n{n-i}\frac{q^{2(i-m)\ceil{\frac{nr}2}}}{1-q^{2(i-m)}}\text{ for }m\geq n.
\end{equation}

All the above formulae are sufficient for calculating the intersection number $\NNN=\intnu$ in the case of (*). The application for $h=2$ case is introduced in Section \S\ref{chetui}. The reader willing to accept these formulae may skip the rest of the section. 

\subsection{Simplification of the intersection formula}
In this subsection, our goal is to prove the formula \eqref{jiaojiao}. For any integer $h$, let $\ek h$, $\ef{2h}$ be the volume of $\EEEEh h$ and $\gggg$ respectively with the standard Haar-measure(see Section \S\ref{hongzha}). It is well-known that
$$
\ef{2h}=\prod_{i=1}^{2h}(1-q^{-i})\qquad \ek h=\prod_{i=1}^h(1-q^{-2i}).
$$

Suppose $\dd g$ is the normalized Haar-measure on $\gggg$, then the standard Haar-measure is $\ef{2h}\dd g$. From now on, we will use standard Haar-measures for the rest of our discussion. The intersection formula is written with the standard Haar-measure $\dd g$ by

\[equation]{\lb{quanzhi}\intnu = \frac{1}{\ep_{K,h}^2}\int_{\gggg}\Big|\Res(\iv\ACj,\iv\ACg)\Big|_F^{-1}\dd g.}
The integrand only depends on $\iv\ACg$.  Let $\SSS=\SS\cap\gggg$. There is an isomorphism
$$
\SSS\cong\gggg/\EEEEo.
$$
We notice that for any $k_1,k_2\in\EEo$, we have
$$
\iv{k_1\ACg k_2}=\iv\ACg.
$$
To simplify notation, we fix the following function throughout the whole section
$$
f(\iv t)=|\mathrm{Res}(\iv t,\iv\ACj)|_F.
$$
Then we can write $\intnu$ as
\[equation]{\lb{diushoujuan}\[split]{
	\frac{1}{\ep_{K,h}^2}\int_{\gggg}f(\iv\ACg)\dd g&=\frac{1}{\ep_{K,h}^2}\int_{\SSS}\int_{\EEEEo}f(\iv{kt})\dd k\dd t\\
	&=\frac{1}{\ep_{K,h}^2}\int_{\SSS}f(\iv{t})\dd t\cdot\int_{\EEEEo}\dd k\\
	&=\frac{1}{\ep_{K,h}}\int_{\SSS}f(\iv{t})\dd t.\\
	}
}
Since we have
$$\SSS=\coprod_{a=0}^h \SSSSSS rha=\coprod_{a=0}^h\coprod_{c=0}^a\SSSSSS rha\cap\SSSSSSS 0hc,$$
we can write the integral \eqref{diushoujuan} as 
$$
\intnu=\frac1{\ep_{K,h}}\int_{\SSS}\ft\dd t=\sum_{a+b+c=h}I(a,b,c)
$$
where $I(a,b,c)$ equals to
$$
I(a,b,c)=\frac1{\ek h}\int_{\SSSSSS rha\cap\SSSSSSS 0h{c}}f(\iv t)\dd t.
$$
For any three polynomials $P_1,P_2,P_3$, by $\mathrm{Res}(P_1,P_2,P_3)$ we mean the product $\mathrm{Res}(P_1,P_2P_3)\mathrm{Res}(P_2,P_3)$. Using the reduction formula in \eqref{maodian}, we can simplify $I(a,b,c)$ to
\[equation]{\lb{huiguo}\[split]{
	&\frac1{\ek h}\int_{\SSSSSS rha\cap\SSSSSSS 0h{c}}f(\iv t)\dd t\\
	=&\frac1{\ek{a+b}\ek{c}}\int_{ \SSShrr {c}0}\int_{\SSSSSS r{a+b}a\cap\SSShr {a+b}0}f(\iv{t_1}\iv{t_2})|\mathrm{Res}(\iv{t_1},\iv{t_2})|_F\dd t_1\dd t_2\dd t_3\\
	=&\frac1{\ek a\ek b\ek c}\int_{ \SSShrr {c}0}\int_{\SSShrr {b}r\cap\SSShr {b}0}\int_{\SSShr ar}f(\iv{t_1}\iv{t_2}\iv{t_3})|\mathrm{Res}(\iv{t_1},\iv{t_2},\iv{t_3})|_F\dd t_1\dd t_2\dd t_3\\
	=&\frac1{\ek a\ek b\ek c}\int_{ \SSShrr {c}0}\int_{\SSShrr {b}r\cap\SSShr {b}0}\int_{\SSShr ar}\frac{|\mathrm{Res}(\iv{t_1},\iv{t_2},\iv{t_3})|_F}{|\mathrm{Res}(\iv\ACj,\iv{t_1}\iv{t_2}\iv{t_3})|_F}\dd t_1\dd t_2\dd t_3.	}
}
Then the equation \eqref{jiaojiao} is obvious if we can prove the following proposition.
\begin{prop}\lb{doremi}
	Suppose $\ACj$ satisfies the condition (*) and $t_1\in\SSShr ar$, $t_2\in\SSShrr br\cap\SSShr b0$, $t_3\in\SSShrr c0$, then
$$
\frac{|\mathrm{Res}(\iv{t_1},\iv{t_2},\iv{t_3})|_F}{|\mathrm{Res}(\iv\ACj,\iv{t_1}\iv{t_2}\iv{t_3})|_F}=\left|{\iv\ACj^{-a}(1)\iv{t_2}^{-b-c}(1)}\right|_F.
$$
\end{prop}
The rest of this subsection is devoted to proving this proposition. 
\begin{defi}
(In this subsection only) we say that an $F$-coefficient polynomial $P_y$ dominates another $F$-coefficient polynomial $P_x$ if any root $\lambda$ of $P_x$ and $\mu$ of $P_y$ have the property
$$
	\vv_F(1-\lambda^{-1})>\vv_F(1-\mu^{-1}).
$$ 
\end{defi}
Our proof will use the following lemma.
\begin{lem}\lb{adadada}
Suppose $P_x,P_y$ are two $F$-coefficient polynomials. Let $a=\deg(P_x)$ and $b=\deg(P_y)$. If $P_y$ dominates $P_x$, then we have
$$
	|\Res(P_x,P_y)|_F=|P_y^a(1)P_x^b(0)|_F.
$$
\end{lem}
\[proof]{Let $\lambda_1,\cdots,\lambda_a$ be roots of $\iv x$, and $\mu_1,\cdots,\mu_b$ roots of $\iv y$. Then by definition, 
$$
	|\mathrm{Res}(\iv x,\iv y)|_F=\prod_{i=1}^a\prod_{j=1}^b\left|\lambda_i-\mu_j\right|_F.
	$$
Since $\iv x(0)=\prod_{i=1}^a\lambda_i$ and $\iv y(0)=\prod_{j=1}^b\mu_j$, we can write
$$
|\mathrm{Res}(\iv x,\iv y)|_F=|\iv x(0)^b\iv y(0)^a|_F\prod_{i=1}^a\prod_{j=1}^b\left|\frac1{\lambda_i}-\frac1{\mu_j}\right|_F.
$$
The assumption $\vv_F(1-\lambda^{-1})>\vv_F(1-\mu^{-1})$ implies
$$
\left|\frac1{\lambda_i}-\frac1{\mu_j}\right|_F=\left|\left(1-\frac1{\lambda_i}\right)-\left(1-\frac1{\mu_j}\right)\right|_F=\left|1-\frac1{\mu_j}\right|_F.
$$
Therefore we have
$$
|\mathrm{Res}(\iv x,\iv y)|_F=|\iv x(0)^b\iv y(0)^a|_F\left|1-\frac1{\mu_j}\right|^a_F=|\iv x(0)^b\iv y(0)^a|_F\left|\frac{\iv y(1)}{\iv y(0)}\right|^a_F.
$$
This lemma follows.
}
\begin{lem}\lb{kuaidian}
	For any $t\in \SSS$, let $\lambda$ and $\mu$ be roots of $\iv t$ and $\iv\ACj$ respectively. If $\iv\ACj$ satisfies (*), then $1-\lambda^{-1}$ and $1-\mu^{-1}$ have different valuation. In other words, either $\iv t$ dominates $\iv\ACj$, or $\iv\ACj$ dominates $\iv t$.
\end{lem}
\begin{proof}
	For any $t\in\SSS$, we only need to show that any $\lambda$ with $\vv_F(1-\lambda^{-1})=\frac dh$ is not a root of $\iv t$ for any odd integer $d$ coprime to $h$. Indeed, if $\iv t$ has such a root and $d$ is coprime to $h$, then all roots $\lambda$ of $\iv t$ have $\vv_F(1-\lambda^{-1})=\frac dh$. On the other hand, $\iv t$ is the characteristic polynomial of $i_t^2$ as an element of $\EE$. Let $\lambda_1,\cdots,\lambda_h$ be the roots of $\iv t$. Then
$$
	\prod_{i=1}^h\left(1-\lambda_i^{-1}\right)=\mathrm{det}_K(1-i_t^{-2}).
$$
	Using \eqref{lemiao} and \eqref{chaojimiao} we know
$$
	1-i_t^{-2}=(i_t^2-1)\circ i_t^{-2}=-e_t^2\circ i_t^{-2} = (e_t\circ i_t^{-1})^2 = t_\#^2.
$$
Then $\vv_F(\mathrm{det}_K(t_\#^2))=d$ is an odd integer. However, we have $t_\#x_0=x_0^\sigma t_\#$. Let $\sigma\in \gggg$ be the linear transformation of Galois conjugation on $K_1^h$. Then $\sigma x_0=x_0^\sigma\sigma$. So $t_\#\sigma$ commutes with $x_0$. Therefore $t_\#\sigma\in\EEo$ and $\vv_K(\det(t_\#)^2)=\vv_K(\det(t_\#\sigma)^2)=2\vv_K(\det(t_\#\sigma))$ is an even number. This is a contradiction. 
\end{proof}

\[lem]{\lb{nimeimei}Let $\iv x$ be an invariant polynomial for an integral double structure. If all roots $\lambda$ of $\iv x$ satisfy $\vv_F(1-\lambda^{-1})>0$, we have $|\iv x(0)|_F=1$. If all roots $\lambda$ of $\iv x$ satisfy $\vv_F(1-\lambda^{-1})\leq 0$, we have $|\iv x(1)|_F=1$.
	}
\begin{proof}Let $\lambda$ be any root of $\iv x$.  If we have $\vv_F(1-\frac1\lambda)>0$, then we must have $|\lambda|_F=1$. Therefore $\iv x(0)$ is the product of all eigenvalues. So $|\iv x(0)|_F=1$.  

If $\vv_F(1-\frac1\lambda)\leq0$, so $\vv_F(\lambda-1)\leq\vv_F(\lambda)$, then by the triangle inequality we have $$\vv_F(1-\lambda)\leq 1$$ for any eigenvalue $\lambda$ of $\iv x$. This implies 
$$
|\iv x(1)|_F\geq 1.
$$Since $\iv x$ is the characteristic polynomial of $i_x^2\in\dEEEo$, the value $\iv x(1)$ is the determinant of
$$
1-i_x^2=e_x^2\in\dEEEo
$$
	by using \eqref{lemiao}. So $|\iv x(1)|_F\leq 1$. This completes the proof.
\end{proof}
Now we are ready to prove Proposition \ref{doremi}.
\begin{proof}[Proof of Proposition \ref{doremi}] Since we have 
$$
t_1\in \SSShr ar,\qquad t_2\in \SSShrr {b}r\cap\SSShr {b}0,\qquad t_3\in \SSShrr {c}0,
$$
$\iv{t_3}$ dominates $\iv{t_2}$. By Lemma \ref{kuaidian}, we have $\iv{t_2}$ dominates $\iv\ACj$. Moreover, $\iv\ACj$ dominates $\iv{t_1}$. By Lemma \ref{adadada}, this implies
	$$|\mathrm{Res}(\iv{t_1},\iv{t_2},\iv{t_3})|_F=|\iv{t_1}(0)^{b+c}\iv{t_2}(0)^c\iv{t_2}(1)^a\iv{t_3}(1)^{a+b}|_F.$$
Similarly,
$$
	|\mathrm{Res}(\iv\ACj,\iv{t_1}\iv{t_2}\iv{t_3})|_F=|\iv\ACj(0)^{b+c}\iv\ACj(1)^a\iv{t_1}(0)^h\iv{t_2}(1)^h\iv{t_3}(1)^h|_F.
$$
	By Lemma \ref{nimeimei}, we have $|\iv{t_1}(0)|_F=|\iv{t_2}(0)|_F=|\iv{\ACj}(0)|_F=|\iv{t_3}(1)|_F=1$. Therefore
$$\frac{|\mathrm{Res}(\iv{t_1},\iv{t_2},\iv{t_3})|_F}{|\mathrm{Res}(\iv\ACj,\iv{t_1}\iv{t_2}\iv{t_3})|_F}=\left|{\iv\ACj^{-a}(1)\iv{t_2}^{-b-c}(1)}\right|_F$$
as desired.\end{proof}
Combining these lemmas, we have obtained a proof of Proposition \ref{doremi}. 
\subsection{Recursion Formula}
Our next goal is to calculate the integrals defined in \eqref{aaa},\eqref{bbb} and \eqref{ccc}. We establish the equation \eqref{a4} in Proposition \ref{yayawawa} and equation \eqref{a3} in Corollary \ref{alala}.
\begin{prop}\lb{yayawawa}
We have
\begin{equation}\lb{girls}
	\sum_{i=0}^a \AAZ i0\BBB{a-i} = \frac{\ep_F}{\ek a^2}=\prod_{i=1}^a\frac{1-q^{1-2i}}{1-q^{-2i}}.
\end{equation}
\end{prop}
\begin{proof}On one hand, by applying our formula \eqref{maodian}, we have
	\[equation*]{\[split]{
		\frac1{\ek a}\int_{\SSSh a}\dd t=&\sum_{i=0}^a\frac1{\ek a\ek{h-a}}\int_{\SSShrr {a-i}0}\int_{\SSShr i0}|\mathrm{Res}(\iv{t'},\iv{t''})|_F\dd t'\dd t''\\
		=&\sum_{i=0}^a\frac1{\ek a\ek{h-a}}\int_{\SSShrr {a-i}0}\int_{\SSShr i0}\left|\iv{t'}(0)^{a-i}\iv{t''}(1)^{i}\right|_F\dd t'\dd t''.
}}
	By Lemma \ref{nimeimei}, we have $\left|\iv{t'}(0)^{a-i}\iv{t''}(1)^{i}\right|_F=1$ for $t'\in\SSShr i0$ and $t''\in\SSShrr {a-i}0$. Then the above integral can be simplified as
	$$
\sum_{i=0}^a\frac1{\ek a\ek{h-a}}\int_{\SSShrr {a-i}0}\int_{\SSShr i0}\dd t'\dd t''=\sum_{i=0}^a \AAZ i0\BBB{a-i}.
	$$
	On the other hand, $\Vol(\SSSh a)=\frac{\Vol(\mathbf{G}_{2a}(\OO F))}{\Vol(\mathbf{G}_a(\OO K))}=\frac{\ep_F}{\ek a}$. Then
	$$
	\frac1{\ek a}\int_{\SSSh a}\dd t=\frac{\ep_F}{\ek a^2}
	$$
as desired.\end{proof}

To compute $\AAA ak$, we consider a more general integral
$$
\AAAA a:= \frac1{\ek a}\int_{\SSShr ar}X^\vf x\dd x.
$$
Besides, we deifne 
$$
\CCCC b:=\frac1{\ek b}\int_{\SSShrr {b}r\cap\SSShr {b}0}|\iv x(1)|_F^{-b}X^\vf x\dd x.
$$
Then we can write
$$
\AAA ai = \left.\AAAA a\right|_{X=q^{2i}}
\qquad
\CCC ai = \left.\CCCC a\right|_{X=q^{2i}}.
$$

\begin{lem}\lb{anata}
We have 
\begin{equation}\lb{zhaolaopo}
	\AAAAA a=\sum_{i=0}^a \AAAA i\CCCCC{a-i}{-2a}.
\end{equation}
\end{lem}
\begin{proof}Clearly
	\[equation*]{\[split]{
		\AAAAA a&=\int_{\SSShr a0}X^{\frac{\vv_F(\iv x(1))}2}\dd x\\
		&=\sum_{i=0}^a\int_{\SSShr {a-i}0\cap\SSShrr{a-i}r}\int_{\SSShr ir}X^{\vf{x'}+\vf{x''}}|\mathrm{Res}(\iv{x'},\iv{x''})|_F\dd x'\dd x''\\
		&=\sum_{i=0}^a\int_{\SSShr {a-i}0\cap\SSShrr{a-i}r}\int_{\SSShr ir}X^{\vf{x'}+\vf{x''}}|\iv{x'}(0)^{a-i}\iv{x''}(1)^i|_F\dd x'\dd x''\\
		&=\sum_{i=0}^a\int_{\SSShr {a-i}0\cap\SSShrr{a-i}r}|\iv{x''}(1)|_F^{i-a}|\iv{x''}(1)|_F^aX^{\vf{x''}}\dd x''\int_{\SSShr ir}X^{\vf{x'}}\dd x'\\
		&=\sum_{i=0}^a\int_{\SSShr {a-i}0\cap\SSShrr{a-i}r}|\iv{x''}(1)|_F^{i-a}(q^{-2a}X)^{\vf{x''}}\dd x''\int_{\SSShr ir}X^{\vf{x'}}\dd x'\\
		&=\sum_{i=0}^a\CCCCC{a-i}{-2a} \AAAA i.\\
}}
This completes the proof.\end{proof}
\begin{cor}\lb{alala}
For all $m\geq n$, we have
$$
	\AAZ n{n-m}=\sum_{i=0}^{n} \CCC im\AAA{n-i}{n-m} 
$$
\end{cor}
\[proof]{By Lemma \ref{anata}, we have $\AAAAA n=\sum_{i=0}^n \AAAA{n-i}\CCCCC{i}{-2n}$. The proof follows by setting $X=q^{2n-2m}$. }

Therefore we proved equations \eqref{a3} and \eqref{a4}.

\subsection{Computation for $\AAA nm$}
In previous steps, we have essentially reduced everything to the computation of $\AAAA i$. In this subsection, we will consider the core $\AAAA i$. We will use the \psc for our computation. This subsection is divided into 3 parts. In the first subsubsection, we introduce the integration over $\SS$ using \psc. In the second subsubsection, we use the \psc to compute $\AAAAA i$. We will see $\AAAAA i$ is a rational function. In the last subsubsection, we give the formula for $\AAAA i$.   
\subsubsection{\PSC}\lb{popo}
The \psc is given by a map 
$$
\maps{\bullet_\#}\OSS{\sss}{x}{\xs:=e_x\circ i_x^{-1}=(x-\zaba)(x-\zaba^\sigma)^{-1}}
$$
where $\OSS$ is the Zariski dense open subset given by $$\OSS=\{x\in\SS:x-\zaba^\sigma\text{ is invertible }\}$$ and $\sss$ is the subspace of $\g$ given by
$$
\sss=\{x_\#\in \g: x_\#\zaba-\zaba^\sigma x_\#=0\}.
$$
We call this map the Polar Stereographic Projection from the pole $\zaba^\sigma$. Note that the set $\sss$ is the set of all semi-$K_1$-linear endomorphisms of $F^{2h}$ with the $K_1$-structure induced by $\baso_1$. 

Our goal for this subsubsection is to prove the following formula.
\begin{prop}
Let $\dd\xs$ be the standard additive Haar-measure for $\sss$ and $\dd x$ the standard Haar-measure for $\SSS$. For any function $\map f\sss{\RR}$, we have
$$
	\int_{\SS} f(\xs)\dd x=\int_{\sss} f(\xs)|\iv x(0)|_F^{2h}\dd \xs.
$$
\end{prop}
\begin{proof}
\newcommand{\tgx}{T_{\ggg,x}}
\newcommand{\tgxz}{T_{\ggg,x_0}}
\newcommand{\tsx}{T_{\SS,x}}
\newcommand{\tsxz}{T_{\SS,x_0}}
\newcommand\tgxs{T_{\g,\xs}}
\newcommand\tsxs{T_{\sss,\xs}}
	
	We consider both $\SS$ and $\sss$ as subvarieties of the ambient varietiy $\g$. For any $x\in\SSS$, let $\tgx$, $\tsx$ be tangent spaces at $x$ for $\g$ and $\SS$. Let $\tgxs$, $\tsxs$ be tangent spaces at $\xs$ for $\g$ and $\sss$. 

	Note that $\dd i_x = \dd e_x = \frac{\dd x}{\zeta-\zeta^\sigma}$, from the following one-form calculation:
\[equation]{\lb{ha}\[split]{
	\dd \xs &=\dd\left(e_x\circ i_x^{-1}\right)\\
	&=e_x\circ i_x^{-1}(\dd i_x) i_x^{-1} + (\dd e_x)i_x^{-1}
	\\&= (e_x\circ i_x^{-1}+1)(\dd i_x)i_x^{-1} \\
	&=(e_x+i_x)\circ i_x^{-2}\circ i_x(\dd i_x)i_x^{-1},
}}
	we deduce that the induced map at the tangent space is given by the following map(we also draw the ambient space in the picture).
$$
	\xymatrix{
\tsx\ar[rr]\ar@{^(->}[d]&&\tsxs\ar@{^(->}[d]\\
		\tgx\ar[rr]&&\tgxs\\
		X\ar@{|->}[rr]&&(e_x+i_x)\circ i_x^{-2}\circ i_x X i_x^{-1}
		}_.
$$

Note that we can factorize this map into the following form

$$
	\xymatrix{
		\tsx\ar[r]\ar@{^(->}[d]&\tsxz\ar[r]\ar@{^(->}[d]&\tsxs\ar@{^(->}[d]\\
		\tgx\ar[r]&\tgxz\ar[r]&\tgxs\\
		X\ar@{|->}[r]&i_x X i_x^{-1}\ar@{|->}[r]&(e_x+i_x)\circ i_x^{-2}\circ i_x X i_x^{-1}\\
		&Y\ar@{|->}[r]&(e_x+i_x)\circ i_x^{-2} Y
		}_.
$$
	Here the first map is induced by the following action of $i_x$
$$
\maps{i_x\cdot}\SS\SS{y}{i_x y i_x^{-1}}
$$
	Since $i_x\in \gg$ and the Haar measure on $\SS$ is $\gg$-invariant, the relative determinant of the first map $X\mapsto i_xXi_x^{-1}$ is 1. This implies the relative determinant for $\mapp\tsx\tsxs$ is equal to the relative determinant of $\mapp\tsxz\tsxs$ with the map $Y\mapsto (e_x+i_x)\circ i_x^{-2} Y$. Note that if we identify $\tgxz$ with $\tgxs$ by additive translation, $\tsxz$ and $\tsxs$ are the same subspaces (Intuitively, the projection plane for the polar-stereographic-projection and the tangent space of the opposite point of the pole are parallel to each other). Since $\tsxz$ is isomorphic to $\Eo$, the relative determinant for $\mapp\tsx\tsxs$ equals
$$
	\mathrm{det}_{K_1}((e_x+i_x)\circ i_x^{-2})^h.
$$
	Note that $(e_x+i_x)^2=1$ implies $|\det_{K}(e_x+i_x)|_K=1$, so
$$
	|\mathrm{det}_{K_1}((e_x+i_x)\circ i_x^{-2})^h|_K=\left|\mathrm{det}_{K_1}\left( i_x^{2}\right)\right|_{K_1}^{-h}=|\iv x(0)|_{K_1}^{-h}=|\iv x(0)|_F^{-2h}.
$$
Now let $\dd x$ and $\dd \xs$ be standard Haar-measures for $\SS$ and $\sss$, we have
$$
|\iv x(0)|_F^{-2h}\dd x = \dd \xs.
$$
This yields the formula that for any function $\map f\sss{\RR}$,
$$
	\int_{\SS} f(\xs)\dd x=\int_{\sss} f(\xs)|\iv x(0)|_F^{2h}\dd \xs.
$$
This completes the proof.\end{proof}

\subsubsection{Computation of $\AAAAA i$}
In this subsubsection, we will prove the following formula.
\begin{thm}\lb{degree}
We have
$$
\AAAAA a=q^{-2a}X\prod_{i=1}^a\frac{1}{1-q^{-2i}X}.
$$
\end{thm}
We briefly introduce our method. Firstly, we split $\sss$ into a disjoint union
$$
\sss=\coprod_{a=0}^h\ssssss0ha
$$
of the following subsets
$$
\ssssss 0ha:=\{x_\#\in\sss: \text{ there is exactly }a \text{ many eigenvalues of }x_\#\text{ with the positive valuation}\}.
$$
Then we will prove the following proposition.
\begin{prop}\lb{dongwu} For any $h,a$, we have
	$$\int_{\ssssss0ha}X^{\vv_K(\xs)}\dd \xs=\AAAAA a.$$
\end{prop}
Suppose this proposition has been proved, then we know 
$$
\int_{\sss}X^{\vv_K(\xs)}\dd \xs=\sum_{i=0}^h\int_{\ssssss 0hi}X^{\vv_K(\xs)}\dd \xs=\sum_{i=0}^h\AAAAA i.
$$
This implies
$$\AAAAA a = \int_{\sssh a}X^{\vv_K(\xs)}\dd \xs-\int_{\sssh{a-1}}X^{\vv_K(\xs)}\dd \xs.$$
Suppose we know the following formula.
\[prop]{\lb{erhuo}We have
$$
\frac1{\ep_{K,a}}\int_{\sssss} X^{\vv_K(\xs)}\dd \xs = \prod_{i=1}^a \frac1{1-q^{-2i}X}.
$$
}
Then we just finish the proof for Theorem \ref{degree} by direct computation
$$
\AAAAA a=\prod_{i=1}^a \frac1{1-q^{-2i}X}-\prod_{i=1}^{a-1} \frac1{1-q^{-2i}X}=q^{-2a}X\prod_{i=1}^a\frac{1}{1-q^{-2i}X}.
$$

We will prove Proposition \ref{erhuo} first, then prove Proposition \ref{dongwu}.
\begin{proof}[Proof of Proposition \ref{erhuo}]Consider the Galois conjugation map (F-linear)
$$
	\maps\sigma{K_1^h}{K_1^h}x{\overline x}_.
$$
	Clearly we have $\sigma\in\sss$ and $\det(\sigma)=(-1)^h$. Then we consider an isomorphism
$$
\maps{l(\sigma)}{\sssss}{\EEEh a}{x}{\sigma\circ x}_.
$$
This implies 
$$
\frac1{\ep_{K,h}}\int_{\sssh a} X^{\vv_K(x)}\dd x=\frac1{\ep_{K,h}}\int_{\EEEh a} X^{\vv_K(x)}\dd x.
$$
	Denote the above integral by $F(a,X)$. Let $\dd g$ be the standard Haar-measure on $\EEh a$. Then we have $\dd x=|g|^a_K\dd g$. We may write $F(a,X)$ into
$$
F(a,X)=\frac1{\ep_{K,h}}\int_{\EEh a\cap\EEEh a}X^{\vv_K(g)}|g|_K^a\dd g.
$$

Let $\Gamma\subset \EEh a$ be the subgroup of upper triangular matrices. By Iwasawa decomposition, we can write $\EEh a=\Gamma\EEEEh a$. Then we can view $\EEh a$ as a homogeneous space with the left $\Gamma\times\EEEEh a$-action given by
$$
\mapi{(\Gamma\times\EEEEh a)\cdot \EEh a}{\EEh a}{(\gamma,k)\cdot g}{\gamma gk^{-1}}_.
$$
The stabilizer of each point is isomorphic to a compact subgroup $\EEEEh a\cap\Gamma$. Then $\EEh a$ has an $\Gamma\times\EEEEh a$-invariant Haar-measure. Since this measure has to be unique, this measure coincide with the Haar-measure of $\EEh a$. We choose the identity matrix $I_h\in\EEh a$ as our base point of $\EEh a$. Then for any $g\in\EEh a$, we can write $g=pI_ht$ for some $p\in\Gamma$ and $t\in \EEEEh a$. Therefore, we have
$$
\int_{\EEh a}\one_{\EEEh a}(g)X^{\vv_K(g)}|g|_K^a\dd g=\frac{\int_{\Gamma\times\EEEEh a}\one_{\EEEh a}(pt)X^{\vv_K(pt)}|pt|_K^a\dd p\dd t}{\int_{\Gamma\cap\EEEEh a}\one_{\EEEh a}(g)X^{\vv_K(g)}|g|_K^a\dd g}.
$$
Note that the integrand $\one_{\EEEh a}(g)X^{\vv_K(g)}|g|_K^a$ is equal to 1 for $g\in\EEEEh a$, we have

$$
\int_{\EEh a}\one_{\EEEh a}(g)X^{\vv_K(g)}|g|_K^a\dd g=\frac{\Vol(\EEEEh a)}{\Vol(\EEEEh a\cap\Gamma)}\int_{\Gamma}\one_{\EEEh a}(p)X^{\vv_K(p)}|p|_K^a\dd p
$$
Since $\Vol(\EEEEh a\cap\Gamma)=\ek 1^a$ and $\Vol(\EEEEh a)=\ek a$, we can write
$$
	F(a,X)=\frac1{\ek1^a}\int_{\Gamma}\one_{\EEEh a}(p) X^{\vv_K(p)}|p|_K^a\dd p.
$$ 
	We remind the reader that $\dd p$ is the standard left-Haar-measure of $P$. Let $\Lambda\subset\Gamma$ be the subgroup of diagonal matrices, $U\subset\Gamma$ the subgroup of unipotent matrices. We take the decomposition $p=\delta u$ such that $\delta\in\Lambda$ and $u\in U$. Note that $\det(u)=1$ for all $u\in U$. We can write
$$F(a,X)=\frac1{\ek1^a}\int_\Lambda\int_U\one[\delta u\in \gl_a(\OO K)]|\delta u|^aX^{\vv_K(\delta u)}\dd u\dd\delta=\frac1{\ek1^a}\int_\Lambda f(\delta)|\delta|^aX^{\vv_K(\delta)}\dd\delta$$
where
$$
f(\delta)=\int_U\one[\delta u\in \gl_a(\OO K)]\dd u.
$$
	Here $\map{\one[\bullet]}{\{\text{True},\text{False}\}}{\{0,1\}}$ is the map such that $\one[\text{True}]=1$ and $\one[\text{False}]=0$. Let $\delta_{11},\cdots,\delta_{aa}\in K$ be diagonal entries of $\delta$, $u_{ij}$ the entry of $u$ in i'th row and j'th column. Then $\delta u\in \gl_a(\OO K)$ is equivalent to that $u_{ij}\in\delta_{ii}^{-1}\OO K$. We thus have
$$
f(\delta)=\prod_{i=1}^a\prod_{j=i+1}^a\int_K\one[u_{ij}\in\delta_{ii}^{-1}\OO K]\dd u_{ij}=\prod_{i=1}^a\prod_{j=i+1}^a|\delta_{ii}|_K^{-1}=\prod_{i=1}^a|\delta_{ii}|_K^{i-a}.
$$
Therefore, 
$$
F(a,X)=\frac1{\ek1^a}\int_{\Lambda\cap\gl_a(\OO K)}\prod_{i=1}^a|\delta_{ii}|_K^{i-a}|\delta_{ii}|_K^aX^{\vv_K(\delta_{ii})}\dd\delta.
$$
This equals
$$
	F(a,X)=\prod_{i=1}^a\frac1{\ek1}\int_{\OO K}(q^{-2i}X)^{\vv_K(\delta_{ii})}\dd\delta_{ii}=\prod_{i=1}^a\frac1{1-q^{-2i}X}
$$
as desired.\end{proof}

Our next goal is to prove Proposition \ref{dongwu}. We will use the method in Section \S\ref{hongzha} for our calculation.

\begin{proof}[Proof of Proposition \ref{dongwu}]
This proof follows the strategy in Section \ref{dasha}. We will adapt the conditions to our situation.

\noindent\textbf{Condition 1:} We choose a subgroup $C$ and set up the $C$-equivariant fibration:
	For any $t\in\ssssss0ha$, we can decompose the characteristic polynomial $P(X)$ of $t$ as $P(X)=P_{0}(X)P_{>0}(X)$ such that $P_{0}(0)\in\OO F^\times$ and $P_{>0}(X)\equiv X^a$ modulo $\pi$. Then let $U_t=\ker(P_{>0}(t))$, which is the maximal invariant subspace such that all eigenvectors of $t$ on $U_t$ has eigenvalue $\lambda$ with $\vv_F(\lambda)>0$. Using this way we have defined a map 
\[equation]{\lb{nvren}
\maps p{\ssssss 0ha}{\Gr}{x}{U_x}_.
}

Let $\sso\subset\sss$ be the subset of invertible matrices. Then $\sso$ is a left-homogeneous space for the group $\EEh a\times\EEh a$ with the action given by

$$
\mapi{(\EEh a\times\EEh a)\cdot\sso}{\sso}{(k_1,k_2)\cdot x}{k_1xk_2^{-1}}_.
$$

Furthermore, we choose our subgroup $C$ as the following 
$$\DEEEEh a:=\{(x,x):x\in\EEEEh a\}.$$ 
Then the surjective map in \eqref{nvren} is a $\DEEEEh a$-equivariant map when we consider $\Gr$ as a $\DEEEEh a$-homogeneous space. 
	
	\noindent\textbf{Condition 2:} We choose a subgroup $P$ such that each fiber is a subset of a $P$-homogeneous space: We denote each fiber $p^{-1}(U)$ by $\ttt$. Then $\ttt$ is a subset of
	$$\tT=\{x\in\sss: xU\subset U\}.$$
Clearly $\tT$ is a homogeneous space of 
	$$P=\PKt\times\PKt$$ where $\PKt$ is the stabilizer of $U_t$.

	\noindent\textbf{Condition 3:} Clearly we have 
	$$\dim(\sso)=\dim(\ttt)+\dim(\Gr).$$

	\noindent\textbf{Notation:} Denote the stabilizer of $t\in\sssss$ by $\TEEEEh a$. We have
$$
	\TEEEEh a:=\{(x,t^{-1}xt):x\in\EEEEh a\}.
$$

Let $\LEE$ and $\LPKK$ be the Lie-algebra of $\EEEEh a$ and $\PKt$ respectively. Let $\LDEEEE,\LTEEEE\subset\LEEEE\times\LEEEE$ be sub Lie-algebras corresponding to $\DEEEEh a,\TEEEEh a\subset\EEEEh a\times\EEEEh a$. Let $\LDPK=\LDEEEE\cap\LPKK\times\LPKK$ and $\LTPK=\LTEEEE\cap\LPKK\times\LPKK$. Let $\LUK=\LEE/\LPKK$, $\LDUK=\LDEEEE/\LDPK$ and $\LTUK=\LTEEEE/\LTPK$. Let $\dd u,\dd u^\Delta$, $\dd u^t$ be corresponding Haar-measures for $\LUK\times\LUK$, $\LDUK$ and $\LTUK$ respectively. 
	
By Theorem \ref{fansile}, we have
$$
\LUK\times\LUK=\LDUK\oplus\LTUK.
$$
Let $J(t)\in F$ be the element such that
$$
\dd u=J(t)\dd u^\Delta\dd u^t.
$$
Theorem \ref{fansile} then implies that we have 
$$
\int_{\sss}\one_{\ssssss0ha}(g)X^{\vv_K(g)}\dd g=\int_{\Gr}\int_{\tT}\one_{\ttt}(t)X^{\vv_K(t)}|J(g)|_F\dd t\dd x.
$$

	Now we will calculate $|J(t)|_F$. Let $t',t''$ be induced linear operators of $t$ on $W$ and $V/W$. We have a natural isomorphism
$$
\LUK\cong\Hom(W,V/W).
$$
Then the isomorphism $\mapp{\LDUK\oplus\LTUK}{\LUK\times\LUK}$ is given by
$$
	\maps{s}{\Hom(W,V/W)\times \Hom(W,V/W)}{\Hom(W,V/W)\times\Hom(W,V/W)}{(x,y)}{(x+y,x+t''^{-1}yt')}_.
$$
Then we use $J(t)$ for the determinant of this map. Note that the valuation of all eigenvalues of $t'$ is larger than 0. And the valuation of all eigenvalues of $t''$ is zero. This implies that all valuations of eigenvalues of $s$ is 0. This implies that $|J(t)|_F=1$. Therefore, we have

$$
\int_{\sss}\one_{\ssssss0ha}(g)X^{\vv_K(g)}\dd g=\int_{\Gr}\int_{\tT}\one_{\ttt}(t)X^{\vv_K(t)}\dd t\dd x.
$$
Since the volume of $\Gr$ is $\frac{\ek h}{\ek a\ek{h-a}}$, we have
$$
\frac1{\ek h}\int_{\sss}\one_{\ssssss0ha}(g)X^{\vv_K(g)}\dd g=\frac1{\ek a\ek{h-a}}\int_{\tT}\one_{\ttt}(t)X^{\vv_K(t)}\dd t.
$$
Furthermore, let $\uuu=\{t\in\EEEE:t|_V=\id, t_{W/V}=\id\}$. After we choose a lifting $\map l{W/U}W$, we can write every element $t\in\ttt$ by $t_1t_2u$ where $u\in\uuu$ and $W$ and $l(W/U)$ are invariant subspaces of $t_1$ and $t_2$, such that $t_1$ acts trivially on $l(W/U)$ and $t_2$ acts trivially on $W$. This implies that $t_1\in\ssssss 0ha$ and $t_2\in \sssssssss h{h-a}$. Therefore we have a decomposition
$$
\ttt=\ssssss 0ha\times\sssssssss h{h-a}\times\uuu.
$$
Then we can write
\[equation*]{\[split]{
	&\frac1{\ek a\ek{h-a}}\int_{\tT}\one_{\ttt}(t)X^{\vv_K(t)}\dd t\\
	=\;&\frac1{\ek a\ek{h-a}}\int_{\ssssss 0ha}\int_{\sssssssss h{h-a}}\int_{\uuu}X^{\vv_K(t_1t_2u)}\dd t_1\dd t_2\dd u\\
	=\;&\frac1{\ek{h-a}}\int_{\sssssssss h{h-a}}\dd t_2\cdot\frac1{\ek a}\int_{\ssssss0ha}X^{\vv_K(t_1)}\dd t_1.\\	
}}

By Proposition \ref{erhuo}, we have
$$
\frac1{\ek{h-a}}\int_{\sssssssss h{h-a}}\dd t_2=\left.\prod_{i=1}^a \frac1{1-q^{-2i}X}\right|_{X=0}=1.
$$
Therefore,
$$
\frac1{\ek h}\int_{\ssssss0ha}X^{\vv_K(g)}\dd g=\frac1{\ek a}\int_{\ssss 0a}X^{\vv_K(t)}\dd t.
$$
Our final goal is to prove this quantity equals to $\AAAAA a$. Indeed, since
$$
\xs^2=e_x\circ i_x^{-1}\circ e_x\circ i_x^{-1}=-e_x^2\circ i_x^{-2}=1-i_x^{-2},
$$
$1-\frac1\lambda$ is an eigenvalue of $\xs^2$ for any root $\lambda$ of $\iv x$. This implies $\xs\in\ssss0a$ if and only if $x\in\SSShr a0$. Furthermore, by Lemma \ref{nimeimei}, $|\iv x(0)|_F=1$. Then
\[equation*]{\[split]{
	\frac1{\ek a}\int_{\ssss 0a}X^{\vv_K(t)}\dd t&=\frac1{\ek a}\int_{\ssss 0a}X^{\vv_K(t)}|\iv x(0)|_F^{2h}\dd t\\
	&=\frac1{\ek a}\int_{\SSShr a0}X^{\vv_K(t)}\dd t\\
	&=\AAAAA a
}}as desired.\end{proof}

\subsubsection{Computation of $\AAAA i$}Now we are able to compute $\AAAA i$ by using complex analysis strategies.	
\begin{lem}\lb{fenshoufenshou}
For any $n,r$, the function $\AAAA n$ is a rational function with poles at $X=q^2,q^4,\cdots,q^{2n}$. Let
$
\aaa nm
$ be the residue of $-q^{-2m}\AAAA n$ at $X=q^{2m}$. Then we have
\begin{equation}\lb{xiangsi}
\AAAA n=\sum_{i=1}^n\frac{\aaa ni(q^{-2i}X)^{\ceil{\frac{nr}2}}}{1-q^{-2i}X}.
\end{equation}
Here $\ceil{r}$ means the smallest integer larger than $r$. In other words, $\ceil{r}=n\iff r\in [n-1,n)$.
\end{lem}
\begin{proof}
	Remember that the degree of a rational function is the order of its pole at the infinity (For example, let $P(x),Q(x)$ be two polynomials. The degree of $\frac{P(x)}{Q(x)}$ is $\deg(P(x))-\deg(Q(x))$). By Lemma \ref{anata}, we have
$$
	\AAAA n=\AAAAA n-\sum_{i=0}^{n-1}\AAAA i\CCCCC{n-i}{-2n}.
$$

	Firstly we claim that the degree of $A(a,r,X)$ is at most $\ceil{\frac{nr}2}-1$. We prove it by induction. When $n=0$, the degree of $\AAAA n=1$ is indeed $0=\ceil{0}-1$. Now we assume the induction hypothesis that the degree of $A(i,r,X)$ is at most $\ceil{\frac{ir}2}-1$ for all $i<n$. By Lemma \ref{degree}, the degree of $\AAAAA n$ is at most 0. For each summand $\AAAA i\CCCCC{n-i}{-2n}$, the degree of $\CCCCC{n-i}{-2n}$ is at most $\lfloor\frac{(n-i)r}2\rfloor$, where the symbol $\lfloor r\rfloor$ means the largest integer no larger than r. By induction hypothesis, the degree of $\AAAA i$ is at most $\ceil{\frac{ir}2}-1$. Therefore, the degree of each summand $\AAAA i\CCCCC{n-i}{-2n}$ is at most
$$
\left\lfloor\frac{(n-i)r}2\right\rfloor+\ceil{\frac{ir}2}-1\leq \ceil{\frac{nr}2}-1.
$$
	This proves our claim. Now let $P(X)$ be a rational function such that
$$
\AAAA n=\sum_{i=1}^n\frac{a(n,r,q^{2i})(q^{-2i}X)^{\ceil{\frac{nr}2}}}{1-q^{-2i}X}+P(X).
$$
	Then $P(X)$ has no poles except the infinity. This implies $P(X)$ is a polynomial. Our claim implies that the degree of $P(X)$ is at most $\ceil{\frac{nr}2}-1$. Since by definition of $\AAAA n$, the coefficient for $X^j$ must be 0 for any $j<\frac{nr}2$, this proves $P(X)=0$, which proves this lemma.
\end{proof}
\begin{cor}
We have 
\begin{equation}\lb{cry}
	q^{2m-2n}\prod_{\substack{i=1\\i\neq m}}^n\frac1{1-q^{2m-2i}}=\sum_{i=m}^n \aaa im\CCC{n-i}{2m-2n}.
\end{equation}
\end{cor}
\begin{proof}
	Using Lemma \ref{fenshoufenshou} and Theorem \ref{degree}, we can write \eqref{zhaolaopo} as
$$
q^{-2n}X\prod_{i=1}^n\frac{1}{1-q^{-2i}X}=\sum_{i=0}^n \sum_{j=1}^i\frac{\aaa ij(q^{-2j}X)^{\ceil{\frac{nr}2}}}{1-q^{-2j}X}\CCCCC{n-i}{-2n}.
$$
Multiplying the above equation by $1-q^{-2m}X$, we can write this equation into
$$
	q^{-2n}X\prod_{i=1,i\neq m}^n\frac{1}{1-q^{-2i}X}=\sum_{i=0}^n\aaa im(q^{-2m}X)^{\ceil{\frac{nr}2}}\CCCCC{n-i}{-2n}+(1-q^{-2m}X)o(X)
$$
	where $o(X)$ refers to a function without poles at $X=q^{2m}$.
This Corollary follows by evaluating this equation at $X=q^{2m}$.
\end{proof}

\section{Computation of the intersection number for $h=2$}\lb{chetui}
In this section, we use our algorithm developed in Section \ref{indu} to compute the arithmetic geometric side of the linear AFL for the case of $h=2$. Our result is listed at the end of this section. Those results are written in a form comparable with the analytic side computed in Section \ref{orba}. 

By Theorem \ref{degree}, we have
\begin{equation}\lb{qusi}
	\aaz n{n-m}=q^{-2m}\prod_{\substack{i=0\\i\neq m}}^{n-1}\frac{1}{1-q^{-2(m-i)}}.
\end{equation}
Furthermore, by \eqref{a1} we have $\aaz {n}{n}=\aaa {n}{n}$ for any $n$. Firstly, we have
$$
\aaa {1}{1}=\aaz {1}{1}=1.
$$
Therefore by the equation \eqref{a2}, we have
$$
\AAA 1{1-m} = \frac{q^{-2m\ceil{\frac r2}}}{1-q^{-2m}}.
$$
Plugging in $m=-1$ and $m=1$, we have
\begin{equation}\lb{chaojidashabi}
\AAA 12=\y{-q^{-2+2\hr}}{}\qquad \AAA 10=\frac{q^{-2\hr}}{1-q^{-2}} \qquad \AAZ10=\frac{q^{-2}}{1-q^{-2}}.
\end{equation}
Note that $\CCC {0}{m}=1$ for any $m$. Using \eqref{a3}, we have
$$
\CCC {1}{m}=\AAZ {1}{1-m}-\AAA {1}{1-m}=\y{q^{-2m}}m-\y{q^{-2m\hr}}m.
$$
Evaluating this expression at $m=0$ and $m=1$, we have
\begin{equation}\lb{b1}
\CCC {1}{1}=\y{q^{-2}-q^{-2\hr}}{},\qquad \CCC {1}{0}=\hr-1,\qquad \CCC {1}{-1}=\frac{q^{2\hr-2}-1}{1-q^{-2}}.
\end{equation}
Continue the same process. By \eqref{a1} we have
$$
\aaa {2}{2}=\aaz {2}{2}=\y{-q^{-2}}{}.
$$
Using \eqref{qusi}, we see
$$
\aaz {2}{1} = \y{q^{-2}}{}.
$$
By \eqref{a1}, we have
$$
\aaa {2}{1}=\aaz {2}{1}-\CCC {1}{1}\aaa {1}{1}=\y{q^{-2\hr}}{}.
$$
Again by applying \eqref{a2}, we have
\begin{equation}\lb{shabi}
\AAA {2}{2-m}=\yy{-q^{-2-2m\gr}}2{2m}+\yy{q^{-2\hr+2(1-m)\gr}}2{2(m-1)}.
\end{equation}
Let $m=2$, we have
\begin{equation}
\AAA {2}{0}=\yy{-q^{-2-4\gr}}2{4}+\Y{q^{-2\hr-2\gr}}{}{2}.
\end{equation}
By Formula \eqref{a3}, we obtain
\begin{equation}\lb{lpd}
\CCC {2}{m}=\AAZ{2}{2-m}-\AAA {2}{2-m}-\CCC {1}{m}\AAA {1}{2-m}.
\end{equation} 
From \eqref{b1} and \eqref{chaojidashabi}, we know $\CCC {1}{0}A[1,2]=\left(1-\hr\right)\y{q^{-2+2\hr}}{}$, therefore moving this term to the left and evaluating \eqref{lpd} at $m=0$, we have
\begin{equation}\lb{onmyway}
\CCC {2}{0}+\left(1-\hr\right)\y{q^{-2+2\hr}}{}=\left.\yy{q^{-2-2m\gr}-q^{-2-2m}}2{2m}\right|_{m=0}+\Y{q^{-2-2\hr+2\gr}-q^{-2}}{}2.
\end{equation}
Applying L'Hospital rule, we have
$$
\CCC {2}{0}=\Y{q^{-2-2\hr+2\gr}-q^{-4}}{}2+\y{\left(\hr-1\right)q^{2\hr-2}-\gr q^{-2}}{}.
$$
Now we calculate $\BBB 1$ and $\BBB 2$. By formula \eqref{a4},
$$
\BBB 1=\frac{1-q^{-1}}{1-q^{-2}}-\AAZ 10=\frac{1-q^{-1}-q^{-2}}{1-q^{-2}}.
$$
By formula \eqref{girls}, we have
$$
\BBB 2=\frac{(1-q^{-1})(1-q^{-3})}{(1-q^{-2})(1-q^{-4})}-\BBB 1\AAZ 10-\AAZ 20
$$
and
$$
\BBB 2=\frac{q^{-3}-q^{-2}}{(1-q^{-2})^{2}}+\frac{q^{-6}-q^{-3}+1-q^{-1}+q^{-4}}{(1-q^{-2})(1-q^{-4})}.
$$
Since $\CCC {1}{0}=\hr-1$, we have
$$
\BBB 1+\CCC {1}{0}=\y{-q^{-1}}{}+\hr
$$
whence
$$
\BBB 2+\CCC {1}{-1}\BBB 1+\CCC {2}{0}
$$
equals

%
%

%
%
%
%
%
%
%
%
%
%
%
%
%
%
%

\[equation*]{
\[split]{
\frac{q^{-2-2\hr+2\gr}}{(1-q^{-2})^2}+\hr\frac{q^{2\hr-2}}{1-q^{-2}}-\frac{q^{2\hr-2}}{1-q^{-2}}-\gr\frac{q^{-2}}{1-q^{-2}}-\frac{q^{-6}}{(1-q^{-2})(1-q^{-4})}\\+
(q^{2\hr-2}-q^{-2})\frac{(1-q^{-1}-q^{-2})}{(1-q^{-2})^2}.
}
}

By our formula \eqref{jiaojiao}, the intersection number equals
$$
N(r)=q^{4r}\AAZ 20+q^{2r}\AAZ 10(\BBB 1+\CCC {1}{0})+(\BBB 2+\CCC {1}{-1}\BBB 1+\CCC {2}{0}).
$$

This equals

\begin{equation*}
\begin{split}
&\frac{q^{-2-2\hr+2\gr}}{(1-q^{-2})^2}+\left(\hr-1\right)\frac{q^{2\hr-2}}{1-q^{-2}}-\gr\frac{q^{-2}}{1-q^{-2}}-\frac{q^{-6}}{(1-q^{-2})(1-q^{-4})}\\&+
(q^{2\hr-2}-q^{-2})\frac{(1-q^{-1}-q^{-2})}{(1-q^{-2})^2}
+\frac{\hr q^{-2\hr+2r}}{1-q^{-2}}-\frac{q^{-1-2\hr+2r}}{(1-q^{-2})^2}\\&+\frac{-q^{-2-4\gr+4r}}{(1-q^{-2})(1-q^{-4})}+\frac{q^{-2\hr-2\gr+4r}}{(1-q^{-2})^2}.
\end{split}
\end{equation*}

Since we have $2r-2\gr=-1$, the intersection formula is simplified to
%
%
%
%
%
%
%

\begin{equation}\lb{tersec}
\begin{split}
N(r)=&\frac{q^{-2-2\hr+2\gr}-q^{-4}}{(1-q^{-2})^2}+\left(\hr-1\right)\frac{q^{2\hr-2}}{1-q^{-2}}-\gr\frac{q^{-2}}{1-q^{-2}}\\+&(q^{2\hr-2}-q^{-2})\frac{(1-q^{-1}-q^{-2})}{(1-q^{-2})^2}+\frac{\hr q^{-2\hr+2r}}{1-q^{-2}}.
\end{split}
\end{equation}
We found $N(\frac12)=1$ and $N(\frac 32)=q+2$. Furthermore, we compute
\begin{equation}
\begin{split}
N(r+2)-N(r)=&\frac{q^{-2\hr+2\gr}}{(1-q^{-2})}+\left(\hr-1\right)q^{2\hr}+\frac{q^{2\hr}}{1-q^{-2}}-2\frac{q^{-2}}{1-q^{-2}}\\+&q^{2\hr}\frac{(1-q^{-1}-q^{-2})}{1-q^{-2}}+\hr q^{-2\hr+2r+2}+\frac{q^{-2\hr+2r+2}}{1-q^{-2}}.
\end{split}
\end{equation}
Note that $2\gr=2r+1$. By simplifying this equation, we may write $N(r+2)-N(r)$ as 
\[equation]{\lb{diffe}
q^{-2\hr+2\gr+1}\left(\frac1{1-q^{-1}}+\hr\right)+q^{2\hr}\left(\hr-\frac{1}{1-q^{-1}}\right)+\frac{2(1-q^{2\hr+2})}{1-q^{2}}.
}

In order to compare this result with orbital integrals, we will rewrite this expression to another form in the rest of the section. Please note that our goal here is preparing a result for comparison with the analytic side rather than simplifying the expression. We do this part of computation only after we know the result of the analytic side. For $a=0$ or $a=1$, we have an identity
\[equation]{\lb{coucou}
(q^{-a+1}+q^a)a+\frac{q^{-a}-q^a}{1-q^{-1}}=0.
}
For any $r\in\frac12\ZZ$, we have
$$
2\hr-\gr=1\qquad \text{or} \qquad 0.
$$
Then let $a=2\hr-\gr$ in the equation \eqref{coucou}, we could write
\[equation*]{
q^{\gr}(q^{-2\hr+\gr+1}+q^{2\hr-\gr})\left(2\hr-\gr\right)
+q^{\gr}\frac{q^{-2\hr+\gr}-q^{2\hr-\gr}}{1-q^{-1}}=0.
}
In other words,
\[equation]{\lb{err}
q^{-2\hr+2\gr+1}\left(2\hr-\gr\right)+q^{2\hr}\left(2\hr-\gr\right)
+q^{\gr}\frac{q^{-2\hr+\gr}-q^{2\hr-\gr}}{1-q^{-1}}=0.
}
Now by computing the difference of the equation \eqref{diffe} and the equation \eqref{err}, we can write $N(r+2)-N(r)$ into

\[equation*]{\[split]{
q^{-2\hr+2\gr+1}\left(\frac1{1-q^{-1}}-\hr+\gr\right)+q^{2\hr}\left(\gr-\hr-\frac{1}{1-q^{-1}}\right)\\+\frac{2(1-q^{2\hr+2})}{1-q^{2}}-q^{\gr}\frac{q^{-2\hr+\gr}-q^{2\hr-\gr}}{1-q^{-1}}.
}}
Simplifying this expression, we obtain 
\[equation]{\lb{shenmewanyi}
N(r+2)-N(r)=2\frac{1-q^{2\hr}}{1-q^2}+\left(\gr+2-\hr\right)q^{2\hr}+\left(\gr-\hr+1\right)q^{2\gr-2\hr+1}.
}
Note that if $2r$ is an odd number, then $\gr=r+\frac12$. If $2r\equiv 1$ mod $4$, we have $\hr=\frac r2+\frac34$. If $2r\equiv 3$ mod $4$, we have $\hr=\frac r2+\frac14$. This implies 
\[equation*]{\left\{\[split]{
	N(r+2)-N(r)=2\frac{1-q^{r+\frac32}}{1-q^2}+\left(\frac r2+\frac74\right)q^{r+\frac32}+\left(\frac r2+\frac34\right)q^{r+\frac12}& \quad \text{ when }2r\equiv 1 \text{ mod }4\\
	N(r+2)-N(r)=2\frac{1-q^{r+\frac12}}{1-q^2}+\left(\frac r2+\frac94\right)q^{r+\frac12}+\left(\frac r2+\frac54\right)q^{r+\frac32}& \quad \text{ when }2r\equiv 3 \text{ mod }4\\
}\right.}

These expressions determine the value of $N(r)$ completely with the initial condition
\[equation*]{\left\{\[split]{
	N\left(\frac12\right)&=1,\\
	N\left(\frac32\right)&=q+2.\\
	}\right.}
%
%

\section{Computation of Orbital Integrals for $h=2$}
\lb{orba}

In this section, we will finish our proof of the linear AFL in $h=2$ case by computing the Analytic side of the linear AFL conjecture. Our test function is the unit of the spherical Hecke algebra. This section have 5 parts. In subsection \S\ref{dafi},\S\ref{lati} and \S\ref{fra}, we briefly describe the combinatorial method for general $h$. Along the way we provide a combinatorial picture of orbital integrals. Subsection \S\ref{divi} describes the orbits that occur in the linear AFL. In subsection \S\ref{compi} we specialize to the computation to the case $h=2$. 

\subsection{Definition of orbital integrals}\lb{dafi}
In this section, we will define the relative orbital integral in \eqref{lashi}. We introduce two important operators in \eqref{centralizer} and \eqref{semilinear}. These two operators will play a central role in our calculation.

Now we prepare materials to define the relative orbital integral with respect to two $\OO F$-algebra embeddings
\[equation*]{\lb{pasa}\[split]{
\map{\blin_1}{\OO F\times \OO F}{\gl_{2h}(\OO F)},\\
\map{\blin_2}{\OO F\times \OO F}{\gl_{2h}(\OO F)}.
}}

Let $g\in\gggg$ be an element such that $\blin_2(x)=g\blin_1(x)g^{-1}$ for any $x\in\OF$. Let $\WWi\subset \gg$ be centralizers of $\blin_i$. It is clear that $\WWi\cong\WW$ for $i=1,2$. For $x\in\WWo$, we write $x=(x_1,x_2)$ for $x_i\in \mathbf{G}_h(F)$. Moreover, define
$$
|x|:=\left|{\det(x_1^{-1}x_2)}\right|_F\qquad \eta_{K/F}(x):=\eta_{K/F}(\det(x_1x_2))
$$
where $\eta_{K/F}$ is the quadratic character of $K/F$. Let $\ooo$ be the characteristic function of $\gggg$. The relative orbital integral is defined by
\[equation]{\lb{lashi}
	\OOOO\blin(\ooo,s):=\int_{{\mathrm C(\blin_1)\cap \mathrm C(\blin_2)}\backslash{\hh\times\hh}}\ooo(u_1^{-1}g u_2)\eta_{E/L}(u_2)\left|u_1u_2\right|^s\dd u_1\dd u_2.
}
where the Haar-measure on $\WWo$ and $\WWo\cap\WWt$ are normalized by $\WWo\cap\gggg$ and $\WWo\cap\WWt\cap\gggg$ respectively. In the case that the invariant polynomial of $(\blin_1,\blin_2)$ has distinct roots, the $F$ algebra
$$
L_{\blin_1,\blin_2}:=\{l\in \g: l\blin_i(x)=\blin_i(x)l \text{ for any }i=1,2\text{ and }x\in\OF\}
$$
is a commutative etale algebra over $F$ and one has $L\cap\gg\cong\WWo\cap\WWt$. 

Now we introduce two important operators. Let 
$$
\zeta = (a,b) \in\OO F^\times \oplus \OO F^\times
$$
be a generator.
Recall Definition \ref{ie}. Let $\cen_{\blin_1,\blin_2}$ be an element in $L_{\blin_1,\blin_2}$ defined by
\[equation]{\lb{centralizer}
\cen_{\blin_1,\blin_2}:=i_{\blin_1,\blin_2}^2\in L_{\blin_1,\blin_2}.
}
Let $\sem_{\blin_1,\blin_2}$ be an element in $\gg$ defined by
\[equation]{\lb{semilinear}
\sem_{\blin_1,\blin_2}:=i_{\blin_1,\blin_2}\circ e_{\blin_1,\blin_2}\in\gg.
}
The trace of $\sem_{\blin_1,\blin_2}$ is zero because $i_{\blin_1,\blin_2}\circ \sem_{\blin_1,\blin_2}\circ i_{\blin_1,\blin_2}^{-1} = -\sem_{\blin_1,\blin_2}$. In this section, we abbreviate those symbols by $\sem$ and $\cen$.
\subsection{Orbital Integral and Lattice Counting}\lb{lati}
In this subsection, we give a combinatorial formula for orbital integrals in terms of lattices. 
A lattice $\lat\subset F^{2h}$ is an $\OO F$-submodule such that $\lat\otimes_{\OO F}F\cong F^{2h}$. Our orbital integral \eqref{lashi} has a natural interpretation of counting lattices. From now, we fix the embedding $\blin_1$ such that 
$$
\blin_1(\zeta)=\mm{aI_h}{}{}{bI_h}.
$$
We fix a lattice
$$
\latz=\OO F^{2h}\subset F^{2h}.
$$

To translate the orbital integral into an object-counting problem, we study the integrand of \eqref{centralizer}, $$\ooo(u_1^{-1}gu_2)\eta_{E/L}(u_2)\left|u_1u_2\right|^s.$$ We will discuss $\ooo(u_1^{-1}gu_2)$ first. Then study $\eta_{E/L}(u_2)\left|u_1u_2\right|^s$.

Note that the integrand of the orbital integral \eqref{lashi} does not vanish only if
$$
\ooo(u_1^{-1}gu_2)\neq 0,
$$
which is equivalent to $gu_2\latz=u_1\latz$. 
Furthermore, it is straightforward to verify $\lat=gu_2\latz=u_1\latz$ if and only if $\lat$ is closed under both actions of $\blin_1(\zeta)$ and $\blin_2(\zeta)$. Finally, since the integrand is considered under the equivalence of $\mathrm C(\blin_1)\cap \mathrm C(\blin_2)=L^\times$, the object-counting process is considering elements of the following subset
$$
\latset := \left\{\lat\subset F^{2h} \text{lattice}, \blin_i(\zeta)\lat=\lat \text{ for }i=1,2\right\}/L^\times.
$$
Here two lattices $\lato,\latt$ are equivalent if and only if $\lato=l\latt$ for some $l\in L^\times$.

Now we explain the term $\eta_{E/L}(u_2)\left|u_1u_2\right|^s$. Let $\lat\subset F^{2h}$ be a lattice with $\zabi\lat=\lat$, there is a direct sum decomposition of $\lat$ corresponding to eigenspaces of $\zabi$ by
$$
\lat = \lat_+\oplus\lat_-
$$
where
\[equation*]{\[split]{
	\lat_+=\im\left(\mapa{\zabi-a}\lat\lat\right)=\ker\left(\mapa{\zabi-b}\lat\lat\right),\\
\lat_-=\im\left(\mapa{\zabi-b}\lat\lat\right)=\ker\left(\mapa{\zabi-a}\lat\lat\right).
}}
The action of $\cen$ preserves components whence $\cen(\lat_+)=\lat_+$ and $\cen(\lat_-)=\lat_-$. We denote their induced maps by $\map{\cen_+}{\lat_+}{\lat_+}$, $\map{\cen_-}{\lat_-}{\lat_+}$. The action of $\sem$ interchanges components. In other words, $\sem(\lat_+)=\lat_-$ and $\sem(\lat_-)=\lat_+$. We denote their induced maps by $\map{\sem_-}{\lat_+}{\lat_-}$, $\map{\sem_+}{\lat_-}{\lat_+}$. Finally, for each representative $\Lambda$ of a lattice class in $\latset$, we can associate to it an order $\sta\lat\subset L$ defined by
$$
\sta\lat=\{l\in L:l\lat\subset \lat\}.
$$
This definition does not depend on the choice of representatives. 

With above constructions, we can write orbital integrals in a combinatorial way as in the following theorem.
\[thm]{\lb{orb}Assume $\latzm=\sem\latzp$, we have
	\[equation]{\lb{quandoushilaji}
	\OOOO\blin(\ooo,s)=\sum_{[\lat]\in \latset}[\OLX:\stax\lat](-q^{ks})^{\mathrm{length}(\lat_+/\sem\lat-)}.
	}
	with $k=2$ if $i_{\tau_1,\tau_2}\in\OO L^\times$ or $k=0$ if $e_{\tau_1,\tau_2}\in\OO L^\times$.
	}
\begin{proof}Consider a subset $\satset\subset\gg/\gggg$ defined by
$$
\satset=\{u\in\gg/\gggg: \zobi u\latz=u\latz \text{ for }i=1,2\}.
$$
	Since $\zabi u_1\latz=u_1\latz$, and $\zibi gu_2\latz=gu_2\latz$ (because $g\zabi g^{-1}=\zibi$), we have 
$$\ooo(u_1^{-1}gu_2)\neq0\implies gu_2\latz=u_1\latz \implies u_1\in\satset.$$
This implies that we can write the orbital integral into
	\[equation*]{
	\int_{\substack{L^\times\backslash\satset\\u_1\latz=gu_2\latz}}|u_1u_2|^s\eta_{K/F}(u_2)\dd u_1\dd u_2=\sum_{\substack{u_1\in\satset\\u_1\latz=gu_2\latz}}[\OO L^\times:\mathrm{Stab}_{L^\times}(u)]|u_1u_2|^s\eta_{K/F}(u_2).
	}
Note that here we can omit $u_1\latz=gu_2$ because $\zaba u_1\latz=u_1\latz$ implies there is an $u_2\in\satset$ with $u_1\latz=gu_2\latz$. Furthermore, this $u_2$ is unique. Then we can write the orbital integral directly into
	\[equation]{\lb{zhongjianshabi}\sum_{\substack{u_1\in\latset}}[\OO L^\times:\mathrm{Stab}_{L^\times}(u)]|u_1u_2|^s\eta_{K/F}(u_2)\quad\text{ with }gu_2\latz=u_1\latz.}
Let $\lat=u_1\latz$. Next we will prove 
 \[equation]{\lb{chara}|u_1u_2|^s\eta_{K/F}(u_2)=(-q^{ks})^{\mathrm{length}(\lat_+/\sem\lat-)}
 }
with $k=2$ if $\irt\in\OO L^\times$ or $k=0$ if $\ert\in\OO L^\times$. For $i=1,2$, denote $\lati=u_i\latz$. Decompose $\lati=\latip\oplus\latim$ by eigenspaces of $\blin_i(\zeta)$. Note that 
	\[equation]{\lb{miejueba}
|u_i|=|\det(u_{i+}^{-1}u_{i-})|_F=\frac{[\latim:\latzm]}{[\latip:\latzp]}
=\frac{[\latim:\sem\latzp]}{[\sem\latip:\sem\latzp]}=[\latim:\sem\latip].
	}
Since $\zibi,\zabi\in\ggg$, we have $i_{\tau_1,\tau_2},e_{\tau_1,\tau_2} \in\ggg$. 

On one hand, If $\det(i_{\tau_1,\tau_2})\in\OO F^\times$, we have $i_{\tau_1,\tau_2}\in\gggg$. Since
$$i_{\tau_1,\tau_2}\zabi=\zibi i_{\tau_1,\tau_2},$$
we have $i_{\tau_1,\tau_2}\latop=\lattp$ and $i_{\tau_1,\tau_2}\lattp=\latop$.  Therefore 
$$
[\lattm:\sem\lattp]=[i_{\tau_1,\tau_2}\lattm:i_{\tau_1,\tau_2}\sem\lattp]=[\latom:\sem\latop].
$$
By \eqref{miejueba}, this implies $|u_1|=|u_2|=q^{\length(\latm:\sem\latp)}$. Furthermore, since 
\[equation]{\lb{quxiao}\eta_{K/F}(u_2)=(-1)^{\vv_F(\det(u_{2+}^{-1}u_{2-}))}=(-1)^{\length(\lattm:\sem\lattp)},}
we have $|u_1u_2|^s\eta_{K/F}(u_2)=(-q^{2s})^{\length(\latm:\sem\latp)}$.

On the other hand, if	$\det(e_{\tau_1,\tau_2})\in\OO F^\times$, we have $e_{\tau_1,\tau_2}\in\gggg$. Since
$$
e_{\tau_1,\tau_2}\zabi=\zibia e_{\tau_1,\tau_2},
$$
we have $\ert\latop=\lattm$ and $\ert\lattm=\latop$.
This implies $|u_2|^{-1}=|u_1|=q^{\length(\latm:\sem\latp)}$. Furthermore, by \eqref{quxiao} again,
we conclude that $|u_1u_2|^s\eta_{K/F}(u_2)=(-1)^{\length(\latm:\sem\latp)}$.

Now after \eqref{chara} has been proved, we could write the orbital integral as

$$
\sum_{\substack{u\in\latset}}[\OO L^\times:\mathrm{Stab}_{L^\times}(u)](-q^{ks})^{\mathrm{length}(\lat_+/\sem\lat-)}.$$
Since $\mathrm{Stab}_{L^\times}(u)$ is the stabilizer of $u$ as an element of $\satset$ under the action of $L^\times$, it is identified with the stabilizer of $\lat=u\latz$ under actions of $L^\times$. Then $\mathrm{Stab}_{L^\times}(u)=\sta\lat^\times$
. This completes the proof of the theorem.
\end{proof}

\subsection{Lattices and Fractional Ideals}\lb{fra}
In this subsection, we simplify the orbital integral in \eqref{quandoushilaji} by parametrizing elements of the set $\latset$ in more details. For each class $[\lat]\in\latset$, each representative $\lat\subset F^{2h}$ is a lattice which is closed under actions of $\zabi$ and $\zibi$. Then we can decompose $\lat=\latp\oplus\latm$ according to eigenspaces of $\zabi$. Each of the component $\latp,\latm$ is a lattice in $F^h$. We denote their stabilizers in $L$ as
\[equation*]{\[split]{
\stap\lat = \{l\in L: l\latp\subset\latp\},\\
\stam\lat = \{l\in L: l\latm\subset\latm\}.\\
}}
It is clear that the definition of $\stap\lat$ and $\stam\lat$ only depends on the class $[\lat]$. We also have $\sta\lat=\stap\lat\cap\stam\lat$. 

Note that $\latp\otimes_{\OO F}F\cong F^{2h}$. A vector $\vv\in F^{h}$ gives rise to an isomorphism $F^{h}\cong L$, which identifies $\latp$ as a $\stap\lat$-submodule of $L$. This is also called a $\stap\lat$-fractional ideal in $L$ in the literature. Two lattices $\latop$ and $\lattp$ are related by $\lattp=l\latop$ for an $l\in L$ if and only if $\stap\lato=\stap\latt$ and they correspond to $\stap\lato$-fractional ideals in the same ideal class.  Let $\mathcal C$ be the set of all fractional-ideal-classes of $L$. Let $\mathcal C(R)\subset \mathcal C$ be the set of $R$-fractional-ideal-classes. 

To describe elements in $\latset$, we need more definitions. Let
$$
\sitset=\{\map\sema\latp\latm\}/\cong,
$$
where two maps $\mapp\latop\latom$, $\mapp\lattp\lattm$ are defining the same element in $\sitset$ if and only if there are $\stap\lat$ and $\stam\lat$-isomorphisms $l_+$ and $l_-$ such that the following diagram commutes
$$
\xymatrix{
	\latop\ar[d]_{l_+}\ar[r]^{\sem_{1-}}&\latom\ar[d]^{l_-}\\
\lattp\ar[r]^{\sem_{2-}}&\lattm
\\
}_.
$$
We have the following description of $\latset$.
\begin{lem}\lb{lalanimei}
We have an isomorphism
$$
	\latset\cong\{([\latp],[\latm],\sema):[\latp],[\latm]\in \mathcal C, \sema\in\sitset, \sem^2\latm\subset\sema\latp\}.
$$
\end{lem}
\begin{proof}
	For any $[\lat]\in\latset$, we have $\sema(\latp)\subset\latm$. Since we also have $\seme(\latm)\subset\latp$, we have $\sem^2\latm=\sema\seme(\latm)\subset\sema\latp$. Therefore we could have a map
	$$
	\mapp\latset\{([\latp],[\latm],\sema):[\latp],[\latm]\in \mathcal C, \sema\in\sitset, \sem^2\latm\subset\sema\latp\}.
	$$	
	We need to show this map is well-defined and that this map is an isomorphism. Let $\lato$, $\latt$ be lattices with $[\lato]=[\latt]$. Then there exists $l\in L$ with $\latt=l\lato$. This implies $\latop=l\lattp$ and $\latom=l\lattm$. Then we have $[\latop]=[\lattp]$ and $[\latom]=[\lattm]$. Furthermore, since $l$ commutes with $\sem_{1-}$, and $\sem_{2-}$, we see we defined the same element in $\sitset$. This map is well defined. 

	On the other hand, we can construct the inverse map and prove the inverse map is well defined. For any  $([\latp],[\latm],\sema)$, we define an element $[\lat]\in\latset$ by the following steps. First, we embed $\latp$, $\latm$ into $F^{2h}$ such that the homomorphism $\map\sema\latp\latm$ and $\map\seme\latm\latp$ are induced by $\sem$. Then we define $\lat=\latp\oplus\latm$. 

	To prove this map is well-defined, we need to show that after the embedding, if $\lattp=l_+\latop$, $\lattm=l_-\latom$ and $\sem\circ l_+=l_-\circ\sem$, then $\latt=l\lato$ for some $l\in L$. Indeed, let $l=l_+\oplus l_-$, we have $l$ commutes with $\sem$. Since $l$ preserves the eigenspaces of $\zabi$, then $l$ commutes with $\zabi$. Similarly $l$ commutes with $\zibi$. This implies $l\in L$.
\end{proof}
By this lemma and \eqref{quandoushilaji}, we can write the orbital integral as
\[equation]{\lb{laile}
	\OOOO\blin(\ooo,s)=\sum_{\substack{R_+\ni\cen\\R_-\ni\cen}}\sum_{\substack{\latp\in C(R_+)\\\latm\in C(R_-)}}\sum_{\substack{\sema\in\sitset\\ \sema\latp\supset \sem^2\latm}}[\OLX:(R_+\cap R_-)^\times](-q^{ks})^{\mathrm{length}(\latm/\sema\latp)}.
}

\subsection{Double structures for division algebra}\lb{divi}

Before we calculate the orbital integral, we study the pair $(\blin_1,\blin_2)$ in more details. By our assumption, it matches to a pair in the division algebra $D$ over $F$ given by
\[equation*]{\[split]{
	\map{\biso_1}{\OO {K_1}}{\OO D}\\
	\map{\biso_2}{\OO {K_2}}{\OO D}\\
	}_.}
By our construction, 
\[equation*]{\[split]{
\sem_{\biso_1,\biso_2}=\ero\circ\iro\\
\cen_{\biso_1,\biso_2}=\iro^2\\
	}_.}
In this subsection, we study the properties of those $\sem=\sem_{\biso_1,\biso_2}$ and $\cen=\cen_{\biso_1,\biso_2}$ that arises from double structures of division algebras. 

\[lem]{\lb{commute}Let $D$ be a quaternion algebra over $L$ and $\ACj\in D$. Let $\zeta\in \OO D^\times$ be a trace zero element and $\zeta\gamma=\gamma\zeta^\sigma$. Then $\vv_D(\gamma)$ is an odd integer. }

\[proof]{
	Let $\varpi$ be the uniformizer of $\OO D$. Since $\OO D/\varpi$ is a finite division algebra, it is a field. Therefore the reduction $\overline{\gamma}$ modulo $\varpi$ must commute with all elements in $\OO D/\varpi$, this implies 
	$$\gamma\zeta-\zeta\gamma\in\varpi\OO D$$
	Since $\gamma\zeta^\sigma=\zeta\gamma$, this implies $\gamma\notin\OO D^\times$.  
	Let $\pi$ be an uniformizer in $\OO L$. Since $D$ is a quaternion algebra, we have $\vv_D(\pi)=2$. Suppose $\vv_D(\gamma)$ is even. Hence $\vv_D(\gamma)=\vv_D(\pi^m)$ for some $m$. Then we have $\pi^{-m}\gamma\in\OO D^\times$ and $\pi^{-m}\gamma\zeta=\zeta^\sigma\pi^{-m}\gamma$, contradiction.  
	}
From now to the rest of the paper, we use $\ACj$ to represent the fixed double structure $(\phi_1,\phi_2)$, and use $\ACj_\#$ to represent $e_{\phi_1,\phi_2}\circ i_{\phi_1,\phi_2}^{-1}$.

\[cor]{Let $\gamma_\#=\ero\circ\iro^{-1}=(\biso_1(\zeta)-\biso_2(\zeta))(\biso_1(\zeta)+\biso_2(\zeta))^{-1}$ be the \psc of the double structure arises in a division algebra. Then
$
\vv_L(\gamma_\#^2)
$ is an odd integer.
}
\begin{proof}
It is clear that $\gamma_\#\in D_L$. Since $D_L$ is a quaterenion algebra over $L$, we have
$$
	\vv_L(\gamma_\#^2)=\vv_{D_L}(\gamma_\#).
$$
	Since $\biso_i(\zeta)\in\OO {D_L}^\times$ and $\gamma_\#\biso_i(\zeta)=\biso_i^\sigma(\zeta)\gamma_\#$, by Lemma \ref{commute}, the number $\vv_{D_L}(\gamma_\#)$ must be an odd integer.
\end{proof}
This corollary implies that for any orbit considered in the linear AFL, the number 
$\vv_F(\iv\ACj(1))$
is always an odd integer.

\subsection{Calculation of orbital integral for $h=2$}\lb{compi}

We call a lattice a principal lattice if it corresponds to a principal fractional ideal. In the case of $h=2$, all lattices $\lat\subset F^{2h}$ that are closed under $\zibi$ and $\zabi$ have principal lattices $\latp$ and $\latm$ as its components. For general $h$, given an $\OO F$-order $R\subset L$, the $R$-fractional ideal class group is not necessarily trivial. The following proposition is well-known and we obmit the proof. 

\[prop]{\lb{principal}Let $L/F$ be a quadratic extension of non-archimedean fields. Then any $\OO F$-lattice $\lat\subset L$ is a principal $R$-fractional ideal in $L$ with $R$ of the form $$R=\OO F+\pi^m\OO L$$ for some positive integer $m$. We denote such a ring as $R_m$.}

In the case of $h=2$, Proposition \ref{principal} implies both $\latp$ and $\latm$ are principal lattices, we may furthermore write the orbital integral as
\[equation]{\lb{sanman}
\OOOO\blin(\ooo,s)=\sum_{\substack{R_+\ni\cen\\R_-\ni\cen}}\sum_{\substack{\sema\in\sotset\\\sema R_+\supset \sem^2R_-}}[\OLX:(R_+\cap R_-)^\times](-q^{ks})^{\mathrm{length}(R_-/\sema R_+)}.
}
Here $[R_+]$ and $[R_-]$ are classes of principal $R_+$ and $R_-$ lattices.

\subsubsection{Calculation of orbital integral}

Using the same notation as previous sections, we have $r=\vv_F(\ACj_\#^2)$. Then $2r=\vv_L(\ACj_\#^2)=\vv_{D_L}(\ACj_\#)$. If $\iro\notin\gggg$, by Theorem \ref{orb} we know the orbital integral do not depend on $s$. Then the derivative is zero. In this subsection, we only consider the non-trivial case where $\iro\in\gggg$. Since we have $\ACj_\#\biso_i(\zeta)=\biso_i^\sigma(\zeta)\ACj_\#$, by Lemma \ref{commute} we know $2r$ is an odd integer. Since $\ACj_\#=\ero\circ\iro^{-1}$, $\vv_F(\iro)=0$ and $\ero\circ\iro=\sem$, we have 
$$\vv_F(\sem^2)=\vv_F(\ACj_\#^2)=2r.$$
This corresponds to $k=2$ case in Theorem \ref{orb}. To make our calculations clear we write the orbital integral in the following way
\[equation]{\lb{ashiba}
\OOOO\blin(\ooo,s)=\sum_{\substack{R_+\ni\cen\\R_-\ni\cen}}[\OLX:(R_+\cap R_-)^\times]\ir.
}
with 
\[equation]{\lb{buzhi}
\ir=\sum_{\substack{\sema\in\sotset\\\lat\supset \sem^2R_-}}(-q^{2s})^{\mathrm{length}(R_-/\sema R_+)}.
}

\[lem]{\lb{chaojiwandan}We have a canonical isomorphism 
$$
\sotset\cong R_-/(R_+\cup R_-)^\times
$$
}
\[proof]{$\sotset$ is the set of maps $\map f{R_+}{R_-}$, where two such maps $f_1,f_2$ are equivalent if and only if $g_-\circ f_1\circ g_+=f_2$ for some $R_+$ and $R_-$-isomorphisms $g_+$ and $g_-$. Since all maps are obtained by multiplying by an element, the map $f$ is uniquely determined by $f(1)$. Suppose $f_1$ and $f_2$ are equivalent, we have $g_-(1)f_1(1)g_+(1)=f_2(1)$. Therefore, we have $
\sotset\cong R_-/(R_+\cup R_-)^\times
$.  }

To study the structure of $\sotset$, we consider a valuation map 
$$
\sotset\cong R_-/(R_+\cup R_-)^\times\longrightarrow \OO L/\OO L^\times\cong\ZZ_{\geq 0}.
$$
Let $\ese,\eso$ be the preimage of even and odd integers under this map respectively. In other words,
$$\ese=\{[\sema]\in\sotset: \vv_L(\sema(1)) \text{ is even } \},$$
$$\eso=\{[\sema]\in\sotset: \vv_L(\sema(1)) \text{ is odd } \}.$$
This gives us a partition
$$
\sotset=\ese\coprod\eso.
$$
Let 
\[equation*]{\[split]{
	\cse&=\{(R_+,R_-,\sema):\sema\in\ese\},\\
	\cso&=\{(R_+,R_-,\sema):\sema\in\eso\},
}}
and
\[equation*]{\[split]{
	\csee{\sem^2}=\{(R_+,R_-,\sema)\in\cse: \sem^2R_+\subset \sema R_-\},\\
	\csoo{\sem^2}=\{(R_+,R_-,\sema)\in\cso: \sem^2R_+\subset \sema R_-\}.\\
	}}
Then we can write
$$
\OOOO\blin(\ooo,s)=\sum_{(R_+,R_-,\sema)\in\csee{\sem^2}\coprod\csoo{\sem^2}}[\OLX:(R_+\cap R_-)^\times](-q^{2s})^{\mathrm{length}(R_-/\sema R_+)}.
$$
Since $\vv_L(\sem^2)$ is odd, and any $(R_+,R_-,\sem_-)\in\csee{\sem^2}$ satisfies $\sem^2 R_+\subset\sem_-R_-$, there is an isomorphism $$\mapi{\csee{\sem^2}}{\csoo{\sem^2}}{\sem_-}{\sem_+}$$
	where the map $\map{\sem_+}{R_-}{R_+}$ is defined by the following commutative diagram
	$$
	\xymatrix{
		&R_-\ar[rd]^{\sem_+}&\\
		R_+\ar[ru]^{\sem_-}\ar[rr]^{\sem^2}&&R_+\\
		}
	$$
Furthermore, we noticed that
$$
\length(R_-/\sema R_+)+\length(R_+/\sem_+ R_-)=\length(R_+/\sem^2 R_+)=2r.
$$
Therefore, we have
\[equation*]{\[split]{
	&\sum_{(R_+,R_-,\sem_+)\in\csoo{\sem^2}}[\OLX:(R_+\cap R_-)^\times](-q^{2s})^{\mathrm{length}(R_+/\sem_+ R_-)}
\\=&
\sum_{(R_+,R_-,\sema)\in\csee{\sem^2}}[\OLX:(R_+\cap R_-)^\times](-q^{2s})^{2r-\mathrm{length}(R_-/\sema R_+)}.
}}
Using this expression we can reduce the orbital integral to the following form

\[equation*]{\[split]{
	&\OOOO\blin(\ooo,s)\\=&\sum_{(R_+,R_-,\sema)\in\csee{\sem^2}}[\OLX:(R_+\cap R_-)^\times]\left((-q^{2s})^{2r-\mathrm{length}(R_-/\sema R_+)}+(-q^{2s})^{\mathrm{length}(R_-/\sema R_+)}\right).\\
}}

We can decompose the set $\cse$ by the following
\[lem]{\lb{shijiewandan}
Let 
$$
\esp abc=\{(R_+,R_-,\sema)\in\cse:R_+=\OO F+\pi^a\OO L,R_-=\OO F+\pi^b\OO L,\vv_L(\sema(1))=2(b-c) \}.
$$
We have 
$$
\csee{\sem^2}=\coprod_{\substack{a\geq c;b\geq c\\ a+b-c\leq\frac{2r-1}2}}\esp abc.
$$
}
\[proof]{
We have $\cse = \coprod_{\substack{a,b,c\in\ZZ_{\geq 0}}}\esp abc$. To prove the lemma, we only need to determine the value of $a,b$ and $c$ so that we have some $x\in R_-$ satisfies $\vv_L(x)=2(b-c)$ and
\[equation]{\lb{hawa}
	\OO F+\pi^b\OO L\supset x(\OO F+\pi^a\OO L)\supset \sem^2(\OO F+\pi^b\OO L).
}
Now we will prove the lemma by considering the inclusion $\OO F+\pi^b\OO L\supset x(\OO F+\pi^a\OO L)$ and $x(\OO F+\pi^a\OO L)\supset \sem^2(\OO F+\pi^b\OO L)$ one by one.

Firstly, by considering all possible elements with odd valuations of each subset, it is clear that $2b\leq 2a+(2b-2c)$ is equivalent to the existence of an element $x\in R_-$ with $\vv_L(x)=2(b-c)$ and
$$
\OO F+\pi^b\OO L\supset x(\OO F+\pi^a\OO L).
$$
Secondly, since $\vv_L(\sem^2)=2r$ is odd, we have $x(\OO F+\pi^a\OO L)\subset \sem^2(\OO F+\pi^b\OO L)$ if and only if $(2b-2c)+2a\leq 2r-1$.
This implies the relation \eqref{hawa} holds only for values $a,b,c\in\ZZ_{\geq 0}$ with $a\geq c, b\geq c$, $a+b-c\leq \frac{2r-1}2$.
	}

\[lem]{\lb{word}We have 
$$
\sum_{(R_+,R_-,\sema)\in\esp abc}[\OLX:(R_+\cap R_-)^\times]=q^{a+b-\max\{0,c\}}.
$$
	}
\[proof]{We have $[\OLX:(R_+\cap R_-)^\times]=q^{\max\{a,b\}}$ for all element in $\esp abc$. Futhermore, we note that there is an isomorphism
$$
\mapi{\esp abc }{ R_-^\times\backslash \left(\pi^{a-c}\OLX\cap R_-\right)/R_+^\times}
{([R_+],[R_-],\sema)}{\sema(1)}
$$
where the representative map is $\map{\sema}{R_+}{R_-}$. 
Since the action of $R_+$ and $R_-$ commute, we have either $R_+\supset R_-$ or $R_+\subset R_-$. Hence 
$$
\esp abc \cong \left(\pi^{a-c}\OLX\cap R_-\right)/R_+\cup R_-
$$
The cardinality of $\esp abc$ is the volume of the above set. Then
$$
\# \esp abc = \frac{[\OLX:(R_+\cup R_-)^\times]}{[\pi^{a-c}\OLX:\pi^{a-c}\OLX\cap R_-]}
=
\frac{q^{\min\{a,b\}}}{q^{\min\{b-(b-c),0\}}}.
$$
This completes the proof of the lemma.}

Besides, it is clear that for any $(R_+,R_-,\sema)\in \esp abc$, we have $\mathrm{length}(R_-/\sema R_+)=a+b-2c$. Using Lemma \ref{word}, Lemma \ref{shijiewandan} and Lemma \ref{chaojiwandan}, we can write the orbital integral $\OOOO\blin(\ooo,s)$ into the following form
$$
\sum_{l=0}^{\frac{2r-1}2}q^l\left(
\sum_{c=-\frac{2r-1}2+l}^{-1}\sum_{\substack{a+b=l\\a,b\geq0}}((-q^{2s})^{l-2c}+(-q^{2s})^{2r-l+2c})
+\sum_{c=0}^l\sum_{\substack{a+b=l+c\\a,b\geq c}}((-q^{2s})^{l-c}+(-q^{2s})^{2r-l+c})
\right).
$$
Now by taking the derivative at $s=0$ we can write
\[equation*]{\[split]{
	&\left.\frac{\dd}{\dd s}\right|_{s=0}\OOOO\blin(\ooo,s){}\\=&\sum_{l=0}^{\frac{2r-1}2}q^l\left(
\sum_{c=-\frac{2r-1}2+l}^{-1}\sum_{\substack{a+b=l\\a,b\geq0}}(-1)^{l}(2(l-2c)-2r)
+\sum_{c=0}^l\sum_{\substack{a+b=l+c\\a,b\geq c}}(-1)^{l-c}(2(l-c)-2r)
\right)\\
=&\sum_{l=0}^{\frac{2r-1}2}q^l\left(
\sum_{c=-\frac{2r-1}2+l}^{-1}(-1)^{l}(l+1)(2(l-2c)-2r)
+\sum_{c=0}^l(-1)^{l-c}(l-c+1)(2(l-c)-2r)
\right)\\
=&\sum_{l=0}^{\frac{2r-1}2}q^l\left(
(-1)^{l}(l+1)\left(\frac{2r-1}2-l\right)
+\left( (-1)^l\left(l^2+2l+\frac12-rl-\frac{3r}2\right)-\left(\frac12+\frac r2\right)\right)
\right)\\
=&\sum_{l=0}^{\frac{2r-1}2}\left((-1)^l\left(\frac l2-\frac r2\right)-\left(\frac12+\frac r2\right)\right)q^l.
}}
Use $N'(r)$ to denote the value for the above expression. We have
\[equation*]{\[split]{
N'(r+2)-N'(r)=-\sum_{l=0}^{\frac{2r-1}2}\left((-1)^l+1\right)q^l
+\left((-1)^{r+\frac12}\left(-\frac34\right)-\left(\frac12+\frac {r+2}2\right)\right)q^{r+\frac12}\\
+\left((-1)^{r+\frac32}\left(-\frac14\right)-\left(\frac12+\frac {r+2}2\right)\right)q^{r+\frac32}.
}}

If $2r\equiv 1$ mod 4, $r+\frac12$ odd, and we have 
\[equation*]{\[split]{
	N'(r+2)-N'(r)&=-2\sum_{l=0}^{\frac r2-\frac14}q^{2l}-\left(\frac r2+\frac34\right)q^{r+\frac12}-\left(\frac r2+\frac74\right)q^{r+\frac32}\\
	&=-\left(2\frac{1-q^{r+\frac32}}{1-q^2}+\left(\frac r2+\frac74\right)q^{r+\frac32}+\left(\frac r2+\frac34\right)q^{r+\frac12}\right).
}}

And if $2r\equiv 3$ mod 4, $r+\frac12$ even, we have
\[equation*]{\[split]{
	N'(r+2)-N'(r)&=-2\sum_{l=0}^{\frac r2-\frac34}q^{2l}-\left(\frac r2+\frac94\right)q^{r+\frac12}-\left(\frac r2+\frac54\right)q^{r+\frac32}\\
	&=-\left(2\frac{1-q^{r+\frac12}}{1-q^2}+\left(\frac r2+\frac94\right)q^{r+\frac12}+\left(\frac r2+\frac54\right)q^{r+\frac32}\right).
}}

It also easy to verify $N(\frac12)=1=-N'(\frac12)$ and $N(\frac32)=q+2=-N'(\frac32)$. This proves
$$
N(r)=N'(r).
$$
We finished the proof of the linear AFL of the case $h=2$ for the identity test function.

\bibliography{Bib/bib}

\providecommand{\bysame}{\leavevmode\hbox to3em{\hrulefill}\thinspace}
\providecommand{\MR}{\relax\ifhmode\unskip\space\fi MR }
\providecommand{\MRhref}[2]{%
  \href{http://www.ams.org/mathscinet-getitem?mr=#1}{#2}
}
\providecommand{\href}[2]{#2}
\begin{thebibliography}{Zha17}

\bibitem[Li18]{qirui2017love}
Qirui Li, \emph{An intersection number formula for {CM} cycles in {Lubin-Tate}
  towers}, arXiv preprint arXiv:1803.07553 (2018).

\bibitem[Zha17]{zhang2017twisted}
Wei Zhang, In preparation, 2017.

\end{thebibliography}
\bibliographystyle{amsalpha}
\end{document}